\input amstex
\documentstyle{amsppt}
\NoBlackBoxes
\hcorrection {-0.2cm}
\magnification=1200
\baselineskip=16pt
\font\we=cmb10 at 14.4truept
\topmatter
\title\nofrills {\we $\Omega$-Admissible Theory II:}
\vskip 0.15cm
\centerline{New metrics on determinant 
of cohomology} And \centerline{Their applications to moduli spaces 
of punctured Riemann surfaces} \endtitle
\author {\bf Lin Weng}
\endauthor
\NoRunningHeads
\endtopmatter
\TagsOnRight
\centerline {Department of Mathematics, Graduate School of
Science, Osaka University,}
\centerline {Toyonaka, Osaka 560, Japan}
\vskip 0.30cm
\noindent 
{\bf Abstract.} For singular metrics, Ray and Singer's analytic torsion
formalism  cannot be applied. Hence we do not have the so-called Quillen
metric on determinant of cohomology with respect to a singular metric. In this
paper, we introduce a new metric  on determinant of cohomology by adapting a totally
different approach. More precisely, by strengthening results in the first paper of
this series, we develop an admissible theory for compact Riemann surfaces with 
respect to singular volume forms, with which the arithmetic Deligne-Riemann-Roch
isometry can be established for singular metrics. As an application, we prove the
Mumford type  fundamental relations for metrized determinant line bundles  over
moduli spaces of punctured Riemann surfaces. Moreover, using an idea of
D'Hoker-Phong and Sarnak, we introduce a natural admssible metric associated to
a punctured Riemann surface via the Arakelov-Poincar\'e
volume, a new invariant for a punctured Riemann surface. With this admissible
metric, we make an intensive yet natural study on two K\"ahler forms on the moduli
space of punctured Riemann surfaces associated to the Weil-Petersson metric and the 
Takhtajan-Zograf  metric (defined by using Eisenstein series. Among others,
 we, together with Fujiki, show that the  Takhtajan-Zograf K\"ahler form is indeed
the first Chern form of a certain metrized line bundle). All this finally leads to a
more geometric interpretation of our new determinant metrics in terms of special
values of Selberg zeta functions. We end this paper by proposing an arithmetic
factorization in terms of Weil-Petersson metrics, cuspidal metrics and Selberg zeta
functions, which then serves as the most global picture for viewing Riemann
surfaces.  
\vskip 0.30cm
\noindent 
{\we Contents}
\vskip 0.30cm
\noindent
\S1. Introduction
\vskip 0.30cm
\noindent
\S2. $\omega$-Arakelov metrics and $\omega$-intersection theory
\vskip 0.30cm
\noindent
\S 3. $\omega$-Riemann-Roch metric and its properties
\vskip 0.30cm
\noindent
\S4. $\omega$-Faltings metric
\vskip 0.30cm
\noindent
\S5. New metrics on determinant of cohomology
for singular metrics and Mumford type isometries
\vskip 0.30cm
\noindent
Appendix to \S5. Universal Riemann-Roch Isomorphism
\vskip 0.30cm
\noindent
\S6. Arakelov-Poincar\'e volume and a geometric interpretation of our new
metrics
\vskip 0.30cm
\noindent
\S7. On Takhtajan-Zograf metric over moduli space of punctured Riemann
surfaces
\vskip 0.30cm
\noindent
Appendix: Arithmetic Factorization Theorem in terms of
Intersection
\vskip 0.30cm
\noindent
\S A1.  Degeneration of Weil-Petersson metrics
\vskip 0.30cm
\noindent
\S A2.  Arithmetric Factorization Theorem: a proposal
\vskip 1.0cm
\noindent
{\we \S1. Introduction}
\vskip 0.30cm
\noindent
Over a compact Riemann surface, for any (smooth) Hermitian line bundle, with respect to any 
(smooth) volume form, we may introduce the Quillen metric ([Qu]) on the
corresponding determinant  of cohomology. Essentially, this is because there
exists only discrete spectrum for the  associated Laplacian, so that the
Ray-Singer's zeta function formalism ([RS]) can be applied.  By using Quillen
metrics, we then have the so-called Deligne-Riemann-Roch isometry, or
equivalently, the Riemann-Roch and the Noether isometries ([De2]).
\vskip 0.30cm
On the other hand, we cannot apply the same strategy to compact Riemann surfaces
with respect to singular volume forms, or better, to punctured Riemann surfaces, due 
to the fact that a certain continuous spectrum exists for the corresponding
Laplacian.  Even though, with respect to hyperbolic metrics on Riemann
surfaces of finite volume, along with the same line as compact Riemann
surfaces, we now have the works done by Efrat ([Ef]), Jorgenson-Lundelius
([JL1], [JL2]), and  Takhtajan-Zograf ([TZ1], [TZ2]) on special values  of
Selberg zeta functions, regularized determinants of Laplacians, and Quillen
metrics, previously it remains to be a  very challenging problem to deduce a
general but natural theory from them.
\vskip 0.30cm
Nevertheless, in this paper, we use a quite independent approach to offer a reasonable metric 
theory for punctured Riemann surfaces. Roughly speaking, we take the Riemann-Roch  and
Noether isometries as the motivation and hence as the final goal for developing such
a theory, since we believe that a  good metric theory for punctured Riemann surfaces
should ultimately provide us these two  isometries in a natural way. Put this in a 
more practical term, we go as follows.
\vskip 0.30cm
As stated above, the up-most main difficuty for doing
arithmetic for singular metrics is the
unpleasent presence of the continuous spectrum for
the associated Laplacian.
We solve this  by developing a
general admissible theory with respect to a
possibly singular volume form $\omega$, which stengthens the results in our previous
paper [We1] in an essential way: we not only deal with points at finite places,
where the metric is finite and smooth, we also develop a system to deal with the
cusps, where the metric is singular. (Please see  (2.3.1), (2.3.2), (2.4.1), (2.5.1)
and (2.5.2) for more details.) Similarly as in [We1],
 the key points at this stage are the existence of the so-called $\omega$-Arakelov
metric and various versions of the Mean Value Lemma, which simply claims that even
though we start with totally independent, possibly singular, volume forms, the
corresponding admissible theories are essentially the same. (Please see
Proposition 2.5.1,  Proposition 2.5.3, Proposition 3.3.1 and Corollary 5.1.2 for
more details.)
\vskip 0.30cm
To apply the general admissible theory to
singular hyperbolic metrics, we then encount
with the second main difficulty: there
exists no geometrically natural admissible metric on the canonical
line bundle. Recall that the singular hyperbolic
metric is  natural only when we view it as a
metric on  the logarithmic canonical line
bundle, which consists of the canonical line
bundle and the cuspidal line bundle. By the
obvious reason, the naive metric on the cuspidal
line bundle resulting only the associated
Dirac symbol is useless for our arithmetic and geometric
consideration: from such a naive metric on the cuspidal line bundle,
we cannot get any admissible metric on the canonical line bundle 
via the decomposition of the original
singular hyperbolic metric on the logarithmic
canonical line bundle; while without using admissible metrics on the canonical line
bundle, it is impossible to apply the general admissible theory.
We overcome this by introducing an invariant
called Arakelov-Poincar\'e volume for a (punctured) Riemann
surface, (please see (6.1.8) for more details,) which exposes the deep relation
between the Euclidean aspect (induced
from the associated Jacobian of its smooth
compactification) and the hyperbolic aspect of
such a Riemann surface at their disposal -- Multiplying the Arakelov
metric with respect to the hyperbolic volume form by this invariant, we get a
natural admissible metric on the canonical line bundle, which is
simply the standard hyperbolic metric when the Riemann surface is compact.
(Please see Corollary 6.4.2 and Remark 6.4.1 for more details.)
 In fact, for compact Riemann surfaces,  by using the Mean
Value Lemma and  a result of D'Hoker-Phong [D'HP] and Sarnak [Sa] on  special
values of Selberg zeta functions and regularized determinants of hyperbolic
Laplacians, such an invariant is first introduced in [We1] to measure the
difference between the standard hyperbolic metric and the Arakelov metric with
respect to hyperbolic volume form. 
\vskip 0.30cm
As an application to moduli spaces of punctured Riemann surfaces, we  give 
Mumford type fundamental isometries for determinant line bundles equipped with our
metrics. (Please see Theorem 5.3.1, Theorem 5.4.1 and Theorem 6.3.1 for more
details). As a direct consequence, we show that the
Weil-Petersson K\"ahler form and the Takhtajan-Zograf K\"ahler form (defined by
using Eisenstein series) on the Teichm\"uller space and on the moduli space of
punctured  Riemann surfaces with fixed signature  naturally arise from  our metric
theory.  Indeed, our study explores the true essence of the pioneer work given by
Takhtajan and Zograf [TZ1,2], which motivates our study. Among others, we, together
with Fujiki, show that the  Takhtajan-Zograf K\"ahler form is indeed
the first Chern form of a certain metrized line bundle. (Please see Theorem 7.3.1
for more details.) All this then leads to a more geometric interpretation of our
determinant metric in terms of spectrum theory. (Please see Theorem 7.3.2 and
Theorem 6.4.1 for more details.)
\vskip 0.30cm
Finally, in an appendix, we propose an arithmetic factorization, which is
motivated by the results of Masur [Ma] and Wolpert [Wo2] on Weil-Petersson metrics
over moduli space of compact Riemann surfaces. The key point here is that  the Weil-Petersson
metric and the Takhtajan-Zograf metric are algebraic so they
are naturally associated to line
bundles on the moduli space together with some
smooth metrics, the so-called local potentials.
On the other hand, such line bundles have natural
extensions to the stably compactification of
moduli space of Riemann surfaces in the sense of
Deligne and Mumford ([DM] and [Kn]), so we may
expect that the associated metrics, or clearly,
local potentials, admit continuous extensions to
the boundary too. Thus by noticing that the
Weil-Petersson metric and Takhtajan-Zograf metric are in
the nature of arithmetic intersection again, we
then may further expect that the above
factorization of line bundles and local
potentials to the boundary give us the corresponding line
bundles and local potentials associated to the
Weil-Petersson metric and the Takhtajan-Zograf metric on
the boundary. We anticipate that such a factorization  plays a key role in studying
the global geometry of Riemann surfaces in furture.
\vskip 0.30cm
As for the language, we intentionally use Deligne pairing [De2], which is certainly
a very natural one for our purpose,  despite the fact that such
a formalism is not as popular as  determinant of cohomology. 
\vskip 0.30cm
As a part of my one semester course at Osaka University in 1997/1998, I explained
and refined all the results in this paper with the help of Professor Fujiki.  I
would like to thank Professor Fujiki, Professor Mabuchi  and Professor Miyanishi
for their supports.  Thanks also due to Professor Ueno, 
Professor Kobayashi and Professor Ohsawa for inviting me to speak on the results
of this paper in a series of lectures at the symposium on Arithemetic Geometry and
Painlev\'e Equations, at Nagoya University respectively, due to Professor To for
fruitful discussions at the earlier stage of this research, which were
unfortunately stopped by some evil force.  Finally, I would like to dedicate this
paper to Serge Lang. 
\vskip 0.30cm
\noindent
{\we \S2. $\omega$-Arakelov metrics and $\omega$-intersection theory}
\vskip 0.30cm
\noindent
(2.1) Throughout this paper, we always assume that $M^0$ is a (punctured) Riemann surface of 
genus $q$. Denote its compactification by $M$, and let $M\backslash
M^0=:\{P_1,\dots,P_N\}$.  We will call $P_i$, $i=1,\dots, N$, {\it cusps} of
$M^0$, and $(q,N)$ the {\it signature} of $M^0$.
\vskip 0.30cm
Recall that a Hermitian metric $ds^2$ on $M^0$ is said to be {\it of
 hyperbolic growth near the cusps}, if for each  $P_i, i=1,\dots,N$,
there exists a punctured coordinate disc
$\Delta^*:=\{z\in {\Bbb C}:0<|z|<1\}$ centered at $P_i$ such that for some
constant $C_1>0$,
\par
\noindent
(i) $\displaystyle{ds^2\leq {{C_1|dz|^2}\over {|z|^2(\log |z|)^2}} \ \
\text{on}\ \Delta^*,}$\hfill (2.1.1) 
\par
\noindent
and there exists a local potential function $\phi_i$ on $\Delta^*$ satisfying
$ds^2={{\partial^2\phi_i}\over {\partial z\partial \bar z}}dz\otimes d\bar z$
on $\Delta^*$, and for some constants $C_2,\,C_3>0$,
\par
\noindent
(ii) $\displaystyle{|\phi_i(z)|\leq C_2\text{max}\{1,\log (-\log |z|)\},\
\text{and}}$\hfill (2.1.2)
\par
\noindent
(iii) $\displaystyle{\left|{{\partial \phi_i}\over {\partial z}}\right|,\,
\left|{{\partial \phi_i}\over {\partial\bar z}}\right|\leq {{C_3}\over
{|z|\left|\log|z|\right|}}\ \ \text{on}\ \Delta^*.}$\hfill (2.1.3)
\vskip 0.30cm
\noindent
In this case, we call $ds^2$ a {\it quasi-hyperbolic metric}, which is 
introduced in [TW1]. (See also [Fu], where a general discussion is given.
Indeed, one may view quasi-hyperbolic metrics as the realization of
good metrics of Mumford [Mu1] in dimension one.)
\vskip 0.30cm
For a quasi-hyperbolic metric $ds^2$ over a punctured Riemann surface $M^0$,
it follows easily from (2.1.1) that 
$\text{Vol}(M^0,ds^2)<\infty$. Denote the normalized volume form of $ds^2$ by $\omega$ 
so that $\text{Vol}(M,\omega)=1$. In this paper, $\omega$
always denotes the normalized volume form on $M$ associated to a
smooth metric (on $M$) or associated to a quasi-hyperbolic metric on $M^0$.
\vskip 0.30cm
\noindent
(2.2) In [TW 1, Theorem 1], we show that there exists a unique $\omega$-{\it
Green's function} 
$g_\omega(\cdot,\cdot)$, or {\it the Green's function with respect to
$\omega$}, on
$M^0\times M^0\backslash
\{\text{diagonal}\}$, such that the  following conditions are satisfied:
\par
\noindent
(i) For fixed $P\in M^0$, and $Q\not=P$ near $P$,
$$g_\omega(P,Q)=-\log|f(Q)|^2+\alpha(Q),$$ where $f$ is a local
holomorphic defining function for $P$, and $\alpha$ is some smooth
function defined near $P$;
\par
\noindent
(ii) $\displaystyle{d_Qd_Q^cg_\omega(P,Q)=\omega(Q)-\delta_P;}$
\par
\noindent
(iii) $\displaystyle{ \int_{M}g_\omega(P,Q)\omega(Q)=0;}$ 
\par
\noindent
(iv) $\displaystyle{g_\omega(P,Q)=g_\omega(Q,P)\ \ \text{for}\
P\not=Q;}$
\par
\noindent
(v) $g_\omega(P,Q)$ is smooth on $\displaystyle{M^0\times M^0\backslash
\{\text{diagonal}\}};$
\par
\noindent
(vi) Near each puncture $P_i$ of $M$, $i=1,\dots, N$, there exists a punctured
coordinate neighborhood $\Delta^*$ centered at $P_i$ such that for fixed
$Q\in \Delta^*$, there exists a constant $C>0$ such that $$|g_\omega(Q,z)|\leq
C\,\text{max}\{1,\log(-\log|z|)\}\ \ \ \ \ \text{on}\ \Delta^*.$$
\par
\noindent 
Here $d_Q^c:={{\sqrt{-1}}\over {4\pi}}(\bar\partial_Q-\partial_Q)$ is with respect
to the second variable (so that $d_Qd_Q^c={{\sqrt{-1}}\over
{2\pi}}\partial_Q\bar\partial_Q$),  and $\delta_P$ is the Dirac delta symbol at
$P$. 
\vskip 0.30cm
The proof comes from the following consideration: for the normalized volume
form $\omega$ associated to a quasi-hyperbolic metric $ds^2$  over a
punctured Riemann surface $M^0$, from definition, it is easy to see that
there exists  a unique locally integrable function
$\beta_\omega$ on $M$ such that 
$$dd^c\beta_\omega=\omega-\omega_{\text{can}},\quad\text{and}\quad
\int_M\beta_\omega(\omega+\omega_{\text{can}})=0.\eqno(2.2.1)$$ 
Here $\omega_{\text{can}}$ denotes the canonical volume form on $M$ defined as follows:
denote by $K_M$ the canonical line bundle of $M$. On $H^0(M,K_{M})$, there exists a natural 
pairing $(\phi,\psi)\mapsto {{\sqrt {-1}}\over 2}\int_{M}\phi\wedge\bar\psi$. Fix any 
orthonormal basis $\{\phi_i\}$ of $H^0(M,K_{M})$ with respect to this 
pairing, by definition,
$$\omega_{\text{can}}:={{\sqrt {-1}}\over {2q}}\sum_{j=1}^q\phi_j\wedge
\bar\phi_j.\eqno(2.2.2)$$ Denote by $g(P,Q)$ the Arakelov-Green's function,
i.e., the Green's function with respect to $\omega_{\text{can}}$. Then we
have
\vskip 0.30cm
\noindent
{\bf Lemma 2.2.1} ([TW1]) {\it With the same notation as above, the
function $g_\omega(P,Q)$ defined on $M^0\times M^0\backslash
\{\text{diagonal}\}$ by
$$g_\omega(P,Q)=g(P,Q)+\beta_\omega(P)+\beta_\omega(Q),\eqno(2.2.3)$$
satisfies the above conditions (i)$\sim$(vi).}
\vskip 0.30cm
\noindent
{\it Proof.} One may prove this lemma as in [La2, Chapter II, 
Proposition 1.3]. The full details are given in my Osaka lecture notes [We2].
In fact, we only
need to remark that with the growth conditions of $\beta$ and $d\beta$, the
arguments in the proof of  [La2, Chapter II, Proposition 1.3] involving
Stokes' theorem remain valid by considering small circles of radius $r$
centered at the punctures and then letting $r\to 0$.  
\vskip 0.30cm
\noindent
(2.3) Now we are ready to define the $\omega$-Arakelov metrics on ${\Cal
O}_M(P)$ for any point $P\in M$ and on $K_M$, the canonical line bundle of
$M$. 
\vskip 0.30cm
First of all, for any $P\in M^0$, define a metric $\rho_{\text{Ar};\omega;P}$
on
${\Cal O}_M(P)$ by setting
$$\log\|1_P\|_{\rho_{\text{Ar};\omega;P}}^2(Q):=-g_\omega(P,Q)+\beta_\omega(P)\quad
\text{for}\ Q\not=P\ \text{in}\ M^0.\eqno(2.3.1)$$
 Here $1_P$ denotes the defining section of ${\Cal O}_M(P)$. (Please note in
particular that the constant $\beta_\omega(P)$ is added.)
 Then
$$\eqalign{~&d_Qd_Q^c(-\log\|1_P\|_{\rho_{\text{Ar};\omega;P}}^2(Q))\cr
=&d_Qd_Q^c(g_\omega(P,Q)-\beta_\omega(P))\quad(\text{by}\ (2.3.1))\cr
=&d_Qd_Q^cg_\omega(P,Q)\cr
=&\omega(Q)-\delta_P\quad(\text{by\ 2.2(ii)})\cr
=&\omega(Q)-\delta_{\text{div}(1_P)}.\cr}$$ Hence $c_1({\Cal
O}_M(P),\rho_{\text{Ar};\omega;P})=\omega$. Here $c_1$ denotes the first
Chern form.
\vskip 0.30cm
Secondly, by Lemma (2.2.1) above, we see that
$$-g_\omega(P,Q)+\beta_\omega(P) =-g(P,Q)-\beta_\omega(Q).$$  Thus, {\bf for
any point} {\bf P}$\in${\bf M}, we (may) define a Hermitian  metric
$\rho_{\text{Ar};\omega;P}$ on ${\Cal O}_M(P)$ by setting
$$\log\|1_P\|_{\rho_{\text{Ar};\omega;P}}^2(Q):=-g(P,Q)-\beta_\omega(Q)\quad
\text{for}\ Q\not=P\ \text{in}\ M^0.\eqno(2.3.2)$$ In particular,
this works also for cusps $P_i$, $i=1, \dots, N$. Easily, we see that
$$c_1({\Cal O}_M(P),\rho_{\text{Ar};\omega;P})=\omega.\eqno(2.3.3)$$ We
will call
$\rho_{\text{Ar};\omega;P}$ {\it the $\omega$-Arakelov metric}, or {\it
the
 Arakelov metric with respect to
$\omega$, on} ${\Cal O}_M(P)$.
\vskip 0.30cm
\noindent
(2.4) A Hermitian line bundle $(L,\rho)$ on $M$ is called $\omega$-{\it
admissible}, if
$c_1(L,\rho)=d(L)\cdot
\omega$. Here $d(L)$ denotes the degree of $L$. From (2.3.3), we have the
following
\vskip 0.30cm
\noindent
{\bf Lemma 2.4.1.} {\it With the same notaton as above, 
$({\Cal O}_M(P),\rho_{\text{Ar};\omega;P})$ is
$\omega$-admissible.}
\vskip 0.30cm
Furthermore, by extending $\rho_{\text{Ar};\omega;P}$ linearly on $P$
by using  tensor products, we know that over any line bundle $L$ on $M$,
there exist 
$\omega$-admissible Hermitian metrics, which are parametrized by ${\Bbb R}^+$.
\vskip 0.30cm
For later use, denote $({\Cal O}_M(P),\rho_{\text{Ar};\omega;P})$ by $\underline {{\Cal O}_M(P)}_\omega$,
or simply $\underline {{\Cal O}_M(P)}$ if no confusion arises. If $(L,\rho)$ is an 
$\omega$-admissible Hermitian
line bundle on $M$, we denote $(L,\rho)$  by $\bar L^\omega$ or simply
$\bar L$ by abuse of notation. Similarly, we use $\bar L(\underline P)$
to denote $\bar L\otimes\underline{{\Cal O}_M(P)}$. 
\vskip 0.30cm
Thus, in particular, on the canonical line bundle $K_M$ of $M$, there exist
$\omega$-admissible  Hermitian metrics. But such metrics are far from being
unique. We next make a  certain normalization.
\vskip 0.30cm
On $K_M$, define {\it the $\omega$-Arakelov metric}
$\rho_{\text{Ar};\omega}$, or {\it the Arakelov metric with respect to}
$\omega$ by setting
$$\|h(z)\,dz\|_{\rho_{\text{Ar};\omega}}^2(P):=|h(P)|^2\cdot\lim_{Q\to P}{{|z(P)-z(Q)|^2}\over{e^{-g_\omega(P,Q)}}}\cdot e^{-2q\beta_\omega(P)}\ \text{for}\
P\in M^0.\eqno(2.4.1)$$ Here $h(z)\,dz$ denotes a section of $K_M$. Then we see that
$$\|h(z)\,dz\|_{\rho_{\text{Ar};\omega}}^2(P)=
\|h(z)\,dz\|_{\text{Ar}}^2(P)\cdot e^{(-2q+2)\beta_\omega(P)}.\eqno(2.4.2)$$
Here $\|\,\|_{\text{Ar}}^2$ denotes the (canonical) Arakelov metric on $K_M$. Thus by the fact that
$\|\,\|_{\text{Ar}}^2$ is $\omega_{\text{can}}$-admissible, (see e.g. 
[La2, Chapter IV, Theorem 5.4],) we have
$$\eqalign{~&c_1(K_M,\rho_{\text{Ar};\omega})\cr
=&(2q-2)\omega_{\text{can}}
+dd^c(-[(-2q+2)\beta_\omega])\quad(\text{by}\ (2.4.2))\cr
=&(2q-2)\omega_{\text{can}}+(2q-2)(\omega-\omega_{\text{can}})\quad(\text{by}\
(2.2.1))\cr =&(2q-2)\omega.\cr}$$ So we have the following
\vskip 0.30cm
\noindent
{\bf Proposition 2.4.2.} {\it With the same
notation as above, $(K_M,\rho_{\text{Ar};\omega})$ is
$\omega$-admissible.} 
\vskip 0.30cm
For later use, denote $(K_M,\rho_{\text{Ar};\omega})$ by 
$\underline {K_M}_\omega$, or simply by $\underline {K_M}$
if no confusion arises. Also we denote $(K_M,\rho_{\text{Ar};\omega}\cdot e^{c\over 2})$
(resp. $\underline{K_M}\otimes\underline {\Cal O_M(P)}$) by $\overline
{K_M}^c$ (resp. $\underline {K_M(P)}$) for any constant $c$.
\vskip 0.30cm
We end this subsection by giving a geometric interpretation for the
$\omega$-Arakelov metric
$\rho_{\text{Ar};\omega}$. We begin with a preperation.
\vskip 0.30cm
Let $\bar L$ be an $\omega$-admissible Hermitian line bundle, then for any point $P\in M$, on the restriction 
$L|_P$, we introduce a  metric 
by multiplying the restriction metric from $\bar L$ to $P$ an additional factor
$\text{exp}[d(L)\cdot{1\over 2}\beta_\omega(P)]$, and we will use the symbol
$ \bar L\|_P$ to indicate the vector space $L|_P$ together
with this modification of the metric, and sometimes call it the
$\omega$-{\it restriction} of $\bar L$ at $P$. With this, by using (2.4.2),
(2.2.1), and the fact that the Arakelov metric induces a natural isometry via
the residue map 
$\text{res}:K_M(P)|_P\to {\Bbb C}$, 
we see that the Arakelov metric with respect to
$\omega$ on $K_M$ is the unique metric such that, at each point $P\in M$,
the natural residue map $\text{res}$ induces the
following  {\it $\omega$-adjunction isometry} $$\text{res}:\underline
{K_M(P)}\|_P\to 
\underline{{\Bbb C}}.\eqno(2.4.3)$$ Here $\underline {\Bbb C}$ denotes
the complex plane ${\Bbb C}$  equipped with the ordinary flat metric.
\vskip 0.30cm
\noindent 
(2.5) For any two line bundles $L,\,L'$ on $M$, denote by $\langle
L,L'\rangle$ the  Deligne pairing associated to $L$ and $L'$. In this
subsection, we define an $\omega$-Deligne norm $h_{\text{De},\omega}$ on
$\langle L, L'\rangle$ for any two $\omega$-admissible Hermitian
line bundles $\bar L$ and $\bar L'$.
\vskip 0.30cm
First, let us define the $\omega$-{\it Deligne norm} for $\langle
{{\Cal O}_M(P)},{{\Cal O}_M(Q)}\rangle$ with $P\not= Q\in M^0$,
for $\omega$-Arakelov metrized line bundles $\underline
{{\Cal O}_M(P)}$ and $\underline {{\Cal O}_M(Q)}$,
 by setting
$$\log\|\langle 1_P,1_Q\rangle\|^2_{h_{\text{De},\omega}}:=
-g_\omega(P,Q)+\beta_\omega(P)+\beta_\omega(Q).\eqno(2.5.1)$$
\vskip 0.30cm
Secondly, note that the right hand side of (2.5.1) can be written as 
$-g(P,Q)$, the Arakelov-Green's function for $P$ and $Q$.
Hence, even though (2.5.1) does not make any sense for cusps, but if we
change it to
$$\log\|\langle 1_P,1_Q\rangle\|^2_{h_{\text{De},\omega}}:=-g(P,Q),\eqno(2.5.2)$$ then we have the 
metrized $\omega$-Deligne pairing $\langle\underline {{\Cal
O}_M(P)},\underline {{\Cal O}_M(Q)}\rangle$  for {\bf all} $P\not= Q\in M$.
\vskip 0.30cm
Finally extending $h_{\text{De},\omega}$ by linearity, we  get a
definition for $\omega$-Deligne norm $h_{\text{De},\omega}(\bar L,\bar L')$
on $\langle  L,L'\rangle$ for any two $\omega$-admissible  Hermitian line
bundles $\bar L$ and $\bar L'$ on $M$. By abuse of notation, we denote 
$\Big(\langle L,L'\rangle,h_{\text{De};\omega}(\bar L,\bar L')\Big)$ simply
by
$\langle \bar L,\bar L'\rangle$.
\vskip 0.30cm
\noindent
{\it Remark 2.5.1.} Even though we study the $\omega$-intersection, the 
Arakelov-Green's function is used in an essential way. This is indeed not quite surprising. 
After all, we only define the $\omega$-intersection for {\bf the}  Hermitian
line bundles 
$\underline {{\Cal O}_M(P)}$ and $\underline {{\Cal O}_M(Q)}$ by using $-g(P,Q)$. Put 
this in a more formal manner, we have the following:
\vskip 0.30cm
\noindent
{\bf Proposition 2.5.1.} (Mean Value Lemma I.) {\it For any two normalized
volume forms 
$\omega_1$ and $\omega_2$ on $M$, there exists a natural isometry
$$\langle\underline{{\Cal
O}_M(P)}_{\omega_1},\underline{{\Cal
O}_M(Q)}_{\omega_1}\rangle\simeq \langle\underline{{\Cal
O}_M(P)}_{\omega_2},\underline{{\Cal
O}_M(Q)}_{\omega_2}\rangle\quad\text{for}\ P\not=Q\in M.
\eqno(2.5.3)$$}
\vskip 0.30cm
As a driect consequence of the $\omega$-adjunction isometry (2.4.3), by definition, we have
the following:
\vskip 0.30cm
\noindent
{\bf Proposition 2.5.2.} ($\omega$-Adjunction Isometry) {\it With the same
notation as above, we have the isometry
$$\langle \underline{K_M(P)},\underline{{\Cal
O}_M(P)}\rangle\simeq
\underline {\Bbb C}\quad\text{for\ any}\ P\in M.\eqno(2.5.4)$$}
\vskip 0.20cm
In a similar style, by using (2.2.1) and (2.4.2), we get
\vskip 0.30cm
\noindent
{\bf Proposition 2.5.3.} (Mean Value Lemma II.) {\it With the same notation as above, for any two normalized volume forms 
$\omega_1$ and $\omega_2$ on $M$, there exists a natural isometry
$$\langle \underline{K_M}_{\omega_1},\underline{K_M}_{\omega_1}\rangle 
\simeq 
\langle
\underline{K_M}_{\omega_2},\underline{K_M}_{\omega_2}\rangle.\eqno(2.5.5)$$}
\vskip 0.30cm
\noindent
{\it Proof.} We may assume that $\omega_2$ is simply $\omega_{\text{can}}$. Denote
$\omega_1$ simply by $\omega$. Then
$$\eqalign{\langle&\underline{K_M}_{\omega},\underline{K_M}_{\omega}\rangle\cr
=&\langle\underline{K_M}_{\omega_{\text{can}}}\cdot e^{-(2q-2)\beta_\omega},
\underline{K_M}_{\omega}\rangle\cr
\simeq&\langle\underline{K_M}_{\omega_{\text{can}}},
\underline{K_M}_{\omega}\rangle\cdot e^{-(2q-2)\int\beta_\omega\cdot
c_1(\underline{K_M}_{\omega})}\cr
\simeq&
\langle\underline{K_M}_{\omega_{\text{can}}},
\underline{K_M}_{\omega_{\text{can}}}\cdot e^{-(2q-2)\beta_\omega}\rangle\cdot
e^{-(2q-2)\int\beta_\omega\cdot c_1(\underline{K_M}_{\omega})}\cr
\simeq&
\langle\underline{K_M}_{\omega_{\text{can}}},\underline{K_M}_{\omega_{\text{can}}}\rangle\cdot e^{-(2q-2)\int\beta_\omega\cdot
c_1(\underline{K_M}_{\omega_{\text{can}}})}
\cdot
e^{-(2q-2)^2\int\beta_\omega\cdot \omega}\cr
\simeq&
\langle\underline{K_M}_{\omega_{\text{can}}},\underline{K_M}_{\omega_{\text{can}}}\rangle\cdot e^{-(2q-2)^2\int\beta_\omega\cdot
\omega_{\text{can}}}
\cdot
e^{-(2q-2)^2\int\beta_\omega\cdot \omega}\cr
\simeq&
\langle\underline{K_M}_{\omega_{\text{can}}},\underline{K_M}_{\omega_{\text{can}}}\rangle\cdot
e^{-(2q-2)^2\int\beta_\omega\cdot (\omega_{\text{can}}+\omega)}\cr
\simeq&
\langle\underline{K_M}_{\omega_{\text{can}}},\underline{K_M}_{\omega_{\text{can}}}\rangle\cdot
e^{-(2q-2)^2\cdot 0}\cr
\simeq&
\langle\underline{K_M}_{\omega_{\text{can}}},\underline{K_M}_{\omega_{\text{can}}}\rangle.\cr}$$ This completes the proof.
\vskip 0.20cm
\noindent
{\it Remark 2.5.1.} The above Mean Value Lemma says that even
though we start with totally independent, possibly singular, volume forms, the
corresponding admissible intersections are essentially the same. 
\vskip 0.30cm
As an application to arithmetic surfaces, we see that 
 the  self-intersection of Arakelov canonical divisor can be understood in any
of these
$\omega$-admissible theories. (For the detailed discussion, see e.g. [We1].)
\vskip 0.45cm
\noindent
{\we \S 3. $\omega$-Riemann-Roch metric and its properties}
\vskip 0.30cm
\noindent
(3.1) With the same notation as in \S2, for any line
bundle $L$ on $M$, denote its associated determinant of cohomology, i.e.,
$\text{det} H^0(M,L)\otimes(\text{det} H^1(M,L))^{\otimes -1}$,
 by $\lambda(L)$. Then it is well-known that we have the following canonical
Riemann-Roch isomophism;
$$\lambda(L)^{\otimes 2}\otimes \lambda({\Cal O}_M)^{\otimes -2}\simeq 
\langle L,L\otimes K_M^{\otimes -1}\rangle.\eqno(3.1.1)$$
(See e.g., [De2], or [Ai].)
\vskip 0.30cm
For a fixed normalized volume form $\omega$ on $M$ associated to a quasi-hyperbolic metric, denote by 
$\underline{K_M}$ the $\omega$-Arakelov canonical line bundle $(K_M,\rho_{\text{Ar};\omega})$.
With respect to $\underline {K_M}$, fix a metric $h_0(\underline {K_M})$ on $\lambda({\Cal O}_M)$. 
Then  for any $\omega$-admissible Hermitian line bundle $\bar L$ on $M$, define an $\omega$-determinant metric
$h_{\text{RR};\underline {K_M};h_0(\underline {K_M})}(\bar L)$ on $\lambda(L)$ by the isometry
$$\Big(\lambda(L),h_{\text{RR};\underline {K_M};h_0(\underline {K_M})}(\bar
L)\Big)^{\otimes 2}
\otimes \Big(\lambda({\Cal O}_M),h_0(\underline {K_M})\Big)^{\otimes
-2}:\simeq
\langle \bar L,\bar L
\otimes \underline{K_M}^{\otimes -1}\rangle.\eqno(3.1.2)$$ We call 
$h_{\text{RR};\underline {K_M};h_0(\underline {K_M})}(\bar L)$ on $\lambda(L)$
the $\omega$-{\it Riemann-Roch metric associated to $\bar L$ with respect to
$\underline {K_M}$ and 
$h_0(\underline {K_M})$}. Since for a fixed $\bar L$, with respect to
$\underline {K_M}$ and 
$h_0(\underline {K_M})$, both $\Big(\lambda({\Cal O}_M),h_0(\underline
{K_M})\Big)$ and $\langle \bar L,\bar L
\otimes \underline{K_M}^{\otimes -1}\rangle$ are fixed,  
$h_{\text{RR};\underline {K_M};h_0(\underline {K_M})}(\bar L)$ is well-defined. By abuse of notation,
we denote $\Big(\lambda(L),h_{\text{RR};\underline {K_M};h_0(\underline
{K_M})}(\bar L)\Big)$ simply by $\underline{\lambda(\bar L)}$.
\vskip 0.30cm
The $\omega$-Riemann-Roch metric satisties the following properties, which are very similar to these for Faltings
metrics. (See Theorem 4.1.1 below.)
\vskip 0.30cm
\noindent
{\bf Proposition 3.1.1.} {\it With the same notation as above, we have
\par
\noindent
(F1)  An isometry of $\omega$-admissible Hermitian line bundles $\bar L\to
\bar L'$ induces an isometry from $\underline{\lambda(\bar L)}$ to
$\underline{\lambda(\bar L')}$; 
\vskip 0.10cm
\noindent
(F2)  If the  $\omega$-admissible metric on $L$ is changed by a factor
$\alpha\in {\Bbb R}^+$, then the metric on $\lambda(L)$
 is changed by the factor $\alpha^{\chi(M,L)}$;
\vskip 0.10cm
\noindent (F3) For any point
$P$ on $M$,  put the  $\omega$-Arakelov metric on
${\Cal O}_M(P)$, and take the tensor metric on $L(-P)$. Then the algebraic
isomorphism
$$\lambda(L)\simeq\lambda(L(-P))\otimes L|P$$  induced by the
short exact sequence of coherent sheaves
$$0\to L(-P)\to L\to L|_P\to 0$$ is  an
isometry $$\underline {\lambda(\bar L)}\simeq\underline {\lambda(\bar L\otimes \underline
{{\Cal O}_M(P)}^{\otimes -1})}\otimes  \bar L\|_P.$$
\par\noindent
(F4) (Serre Isometry)
$\displaystyle{\Big(\lambda(K_M),h_{\text{RR};
\underline{K_M};h_0(\underline{K_M})}(\underline{K_M})\Big)
\simeq \Big(\lambda({\Cal O}_M),h_0(\underline{K_M})\Big).}$}
\vskip 0.30cm
\noindent
{\it Proof.}  (F4) is simply the Serre duality. (F1) and (F2) are direct
consequence of  the definitions of the $\omega$-Riemann-Roch metric and 
the $\omega$-intersection. Finally, (F3) is a direct consequence of the definition of the 
$\omega$-Riemann-Roch  metric and the $\omega$-adjunction isometry, which  also explains 
why in our definition of the $\omega$-Riemann-Roch  metric and the proposition here
we use
$\underline{K_M}$ and $\underline{P}$, i.e., $K_M$ and ${\Cal O}_M(P)$ together with
the $\omega$-Arakelov metrics.
\vskip 0.30cm
\noindent
{\it Remark 3.1.1}. By (F4), we see that giving a normalization for
$h_0(\underline{K_M})$ on $\lambda({\Cal O}_M)$ is equivalent to normalizing
$h_{\text{RR};\underline{K_M};h_0(\underline{K_M})}$ on $\lambda(K_M)$.
\vskip 0.30cm
\noindent
(3.2) Similarly, with respect to $\overline {K_M}$, i.e., $K_M$ with
an arbitrary $\omega$-admissible metric, we fix a metric
$h_0(\overline {K_M})$ on $\lambda({\Cal O}_M)$. Then  with respect to $\overline
{K_M}'$, i.e.,
$K_M$ equipped with (possibly) another
$\omega$-admissible Hermitian metric, and $h_0(\overline {K_M})$, for any
$\omega$-admissible Hermitian line bundle $\bar L$, we may define the
associated Riemann-Roch metric, denoted by $h_{\text{RR};\overline
{K_M}';h_0(\overline {K_M})}(\bar L)$, by the isometry
$$\Big(\lambda(L),h_{\text{RR};\overline {K_M}';h_0(\overline {K_M})}(\bar
L)\Big)^{\otimes 2}
\otimes \Big(\lambda({\Cal O}_M),h_0(\overline {K_M})\Big)^{\otimes
-2}:\simeq
\langle \bar L,\bar L
\otimes (\overline{K_M}')^{\otimes -1}\rangle.\eqno(3.2.1)$$
The dependence of $h_{\text{RR};\overline {K_M}';h_0(\overline {K_M})}(\bar L)$ on $\bar L$ and 
$\overline{K_M}'$ is clear, as it is given by the $\omega$-intersection theory. More precisely, 
directly from the defintion, we have
\vskip 0.30cm
\noindent
{\bf Proposition 3.2.1.} {\it The dependence of $h_{\text{RR};\overline {K_M}';h_0(\overline {K_M})}(\bar L)$ on $\bar L$ and 
$\overline{K_M}'$ is given by the following equality:
$$h_{\text{RR};\overline {K_M}'\otimes {\Cal O}_M(e^c);
h_0(\overline {K_M})}(\bar L\otimes {\Cal O}_M(e^f))=
h_{\text{RR};\overline {K_M}';h_0(\overline {K_M})}(\bar L)\cdot 
e^{\chi(L)\cdot f-d(L)c/2}.
\eqno(3.2.2)$$ Here for a constant $c$, ${\Cal O}_M(e^c)$ denotes the trivial
line bundle equipped  with the metric $\|1\|^2=e^c.$}
\vskip 0.30cm
On the other hand, the dependence of $h_{\text{RR};\overline {K_M}';h_0(\overline {K_M})}
(\bar L)$ on $\overline{K_M}$ is not so easy to determined. Indeed, the most essential part
for such a dependence
is independent of the above (weak) Riemann-Roch isometry. Nevertheless,  
from our study on the admissible theory 
with respect to smooth volume forms in [We1],  it is very natural
to take  the following principle, which has its root from the Polyakov
variation formula (see e.g., [Fay2, (3.30)]):
$$h_0(\overline{K_M}^c):=h_0(\underline{K_M})\cdot
e^{{{2q-2}\over{12}}\cdot c}.\eqno(3.2.3)$$ Here, as before,
$\overline{K_M}^c=\underline{K_M}\otimes {\Cal O}_M(e^c).$ 
\vskip 0.30cm
\noindent
{\it Remark 3.2.1.} The reader may ask why we use $2q-2$ instead of using
$2q-2+N$ in (3.2.3). We justify our choice by the following observation: the
Hermitian metric on $K_M$ used in (3.2.3) would have the first Chern form
$(2q-2)\omega$, which is different from the singular metric introduced by the
quasi-hyperbolic metric, for which the total volume is $2\pi(2q-2+N)$. So
our normalization is in the same spirit as the one in [JL2, \S7].
\vskip 0.30cm
That is, we have the following
\vskip 0.30cm
\noindent
{\bf Proposition-Definition 3.2.2.} (Polyakov Variation Formula I) {\it With
the same notation as above,  we have the following relation
$$h_{\text{RR};\overline {K_M}';
h_0(\overline {K_M}\otimes {\Cal O}_M(e^c))}(\bar L)=
h_{\text{RR};\overline {K_M}';h_0(\overline {K_M})}(\bar L)
\cdot e^{{{2q-2}\over {12}}\cdot c}.\eqno(3.2.4)$$}
\vskip 0.20cm
We end this section by the following  consequence of the definition of
the $\omega$-Riemann-Roch metric.
\vskip 0.30cm
\noindent
{\bf Proposition 3.2.3.} (Serre Isometry) {\it With the same notation as above, we get the isometry:
$$\Big(\lambda(L),h_{\text{RR};\overline {K_M}';h_0(\overline {K_M})}(\bar L)
\Big)\simeq \Big(\lambda(K_M\otimes L^{\otimes -1}), 
h_{\text{RR};\overline {K_M}';h_0(\overline {K_M})}
(\overline{K_M}'\otimes \bar L^{\otimes -1})\Big).\eqno(3.2.5)$$}
\vskip 0.20cm 
\noindent
(3.3) In (3.1) and (3.2), for a {\bf fixed} normalized volume form $\omega$
on
$M$, we introduce 
$h_{\text{RR};\overline {K_M}';h_0(\overline {K_M})}(\bar L)$ in such a way that if one of 
$h_0(\overline{K_M}'')$ is fixed, then all other
determinant metrics $h_{\text{RR};\overline {K_M}'; h_0(\overline
{K_M})}(\bar L)$ are fixed, by using (3.2.2) and (3.2.4), or better
Proposition 3.2.1 and Proposition 3.2.2.
\vskip 0.30cm
Now we explain how the $\omega$-Riemann-Roch metrics depend on $\omega$. 
Similarly, motivated by our work  
on admissible theory with respect to smooth volume forms in [We1], we relate
different $\omega$-Riemann-Roch metrics by using the following isometry: for
any two normalized volume forms $\omega_1$ and $\omega_2$ on $M$,
$$\Big(\lambda({\Cal O}_M),h_0(\underline{K_M}_{\omega_1})\Big)\simeq
\Big(\lambda({\Cal O}_M),
h_0(\underline{K_M}_{\omega_2})\Big).\eqno(3.3.1)$$ In other words, {\it even though
$\underline{K_M}_\omega$, the $\omega$-Arakelov canonical line bundle, depends on
$\omega$ in an essential way, but the induced metric on the determinant of
cohomology does not depend on
$\omega$ at all.} We may say that this is one of the most important discovery in
[We1], where we establish this relation for Quillen metrics. 
 As a direct
consequence, we get the following;
\vskip 0.30cm
\noindent
{\bf Proposition 3.3.1.} (Mean Value Lemma III) {\it With the same notation
and normalization as above, for any two normalized volume forms $\omega_1$
and $\omega_2$ on $M$, we get the following isometries:
\vskip 0.30cm
\noindent
(a) (Polyakov Variation Formula II) $$
\underline{\lambda(\underline{K_M}_{\omega_1})}_{\omega_1}
\simeq
\underline{\lambda(\underline{K_M}_{\omega_2})}_{\omega_2}.
\eqno(3.3.2)$$
\vskip 0.30cm
\noindent
(b) For all $n_j\in {\Bbb Z}$ and $Q_j\in M$,
$$\underline{\lambda(
\underline{{\Cal O}_M(\Sigma_jn_jQ_j)}_{\omega_1})}_{\omega_1}
\simeq \underline{\lambda(
\underline{{\Cal
O}_M(\Sigma_jn_jQ_j)}_{\omega_2})}_{\omega_2}.\eqno(3.3.3)$$}
\vskip 0.30cm
\noindent
{\it Proof.} The key point here is that on $K_M$ and ${\Cal O}_M(Q_j)$, all metrics
are carefully chosen to be $\omega$-Arakelov metrics. Thus
 (a) comes from the Serre isometry; while (b)
is deduced from the Riemann-Roch isometry and  the Mean Value Lemma I for
$\omega$-arithmetic intersection  by a tedious calculation.
\vskip 0.30cm
Thus by the above three kinds of normalizations for the $\omega$-determinant
metrics, i.e., (3.2.2), (3.2.4) and (3.3.2), we see that in order to get the
$\omega$-Riemann-Roch metric uniquely, we have two different ways to
normalize them:  one is to uniformly define metrics
$h_0(\underline{K_M}_\omega)$ for all normalized volume forms $\omega$
first, which satisfy (3.3.2), i.e., Proposition 3.3.1(a), and then use
Proposition 3.2.1 and Proposition 3.2.2 to get all other metrics for any
admissible line bundles; while the other is to define for all
$\omega$-admissible Hermitian line bundle $\bar L$ the Riemann-Roch metrics
$h_{\text{RR};\overline{K_M}^\omega;h_0({\overline{K_M}'}^\omega)}(\bar L)$
on $\lambda(L)$ for a certain fixed
$\omega$, then to check these metrics satisfy (3.3.2) and (3.2.4) and hence
is compactible with our theory. We next give two independent approaches to show how
this can be possibly done in a very concrete manner. The first  is with
respect to any normalized volume form, which then gives an alternative way to
find Quillen metric; while the second   works only for singular hyperbolic
volume forms, which then leads to a more geometric interpretation of our new
determinant metric. 
\vskip 0.30cm
\noindent
{\we \S4. $\omega$-Faltings metric}
\vskip 0.30cm
\noindent
(4.1) This approach begins with the following condition.
\vskip 0.30cm
\noindent
{\it (F0) With respect to the normalized volume $\omega$ associated to a
quasi-hyperbolic metric $d\mu$  on a compact Riemann surface $M$, the metric
$h_{\text{RR};\underline{K_M};h_0(\underline{K_M})}$ on $\lambda(K_M)$ is
defined to be the determinant of the Hermitian metric on 
$H^0(M,K_M)$ induced from the following natural pairing $$(\phi,\psi)\mapsto
{{\sqrt{-1}}\over 2}\int_M\phi\wedge\bar\psi.\eqno(4.1.1)$$}
\vskip 0.30cm
\noindent
{\it Remark 4.1.1.} It appears that (F0) is quite stange as no
$\omega$ is involved (in the natural pairing). But one should not understand
in this way, as it is obvious that  the above
natural paring on
$H^0(M,K_M)$ can also be defined by using any metric $d\mu$ on the Riemann
surface $M$ due to the fact that the dimension of the base manifold $M$ is one (so
that the dual of the tangent bundle is simply the canonical line bundle).
\vskip 0.30cm
Now we may improve Proposition 3.1.1 as follows.
\vskip 0.30cm
\noindent
{\bf Theorem 4.1.1.} {\it With respect to the normalized volume $\omega$ 
on a compact Riemann surface $M$, for any
$\omega$-admissible Hermitian  line bundle $\bar L$,
there exists a unique metric $h_{\text{RR};\underline{K_M};h_0(\underline{K_M})}(\bar L)$, denoted also by
$h_{F;\omega}(\bar L)$ and called the $\omega$-Faltings metric, on
$\lambda(L)$ such that conditions (F0) $\sim$ (F5) are satisfied. Moreover,
we have the following Riemann-Roch isometry:
$$\Big(\lambda(L),h_{F;\omega}(\bar L)\Big)^{\otimes 2}
\otimes \Big(\lambda({\Cal O}_M),h_{F;\omega}(\underline {{\Cal
O}_M})\Big)^{\otimes -2}\simeq
\langle \bar L,\bar L
\otimes \underline{K_M}^{\otimes -1}\rangle.\eqno(4.1.2)$$}
\par
\noindent
{\it Proof.} The proof of this theorem is the same as the one for 
Faltings' original theorem [Fa, Theorem 1]. Namely, fixed
a large enough positive integer $r$ and a degree $r+q-1$ divisor $E$ on $M$. Then, for any point
$(Q_1,\dots,Q_r)\in M^r$, by a tedious calculation,
we get the isometry $$\eqalign{~&\Big(\lambda({\Cal O}_M(E-Q_1-\dots-Q_r)),
h_{F;\omega}(\underline{{\Cal O}_M(E-Q_1-\dots-Q_r)})\Big)\cr
\simeq&\Big(\lambda({\Cal O}_M(E)),h_{F;\omega}(\underline{{\Cal
O}_M(E)})\Big)\otimes
\otimes_{i=1}^r\Big({\underline{{\Cal O}_M(E)}}\|_{Q_i}\Big)^{\otimes -1}\cr
&\otimes
\otimes_{1\leq i<j\leq r}\Big({\underline{{\Cal
O}_M(Q_i)}}\|_{Q_j}\Big),\cr}$$ by using (F3). Noticing that in Faltings' original
theorem, when dealing with the norm on the restriction, Faltings uses the most
direct restriction for the canonical volume form, while for us, we modify it to
$\|$, i.e., we use the
$\omega$-restriction. Thus after taking $dd^c$ on
$M^r$,  $$c_1({\underline{{\Cal O}_M(E)}}\|_{Q_i})=(r+q-1)
\text{pr}_i^*(\omega-dd^c\beta_\omega)=(r+q-1)\text{pr}_i^*\omega_{\text{can}}.$$
Here
$\text{pr}_i:M^r\to M$ denotes the $i$-th projection. Similarly,
$$c_1({\underline{{\Cal
O}_M(Q_i)}}\|_{Q_j})=dd^c\Big(g_\omega(Q_i,Q_j)-\beta_\omega(Q_i)\Big)
-dd^c\beta_\omega(Q_j)=dd^cg(Q_i,Q_j).$$ Thus, we get
$$\eqalign{~&c_1\Big(\lambda({\Cal
O}_M(E-Q_1-\dots-Q_r)),h_{F;\omega}(\underline{{\Cal
O}_M(E-Q_1-\dots-Q_r)})\Big)\cr
=&-(r+q-1)\sum_{i=1}^r\text{pr}_i^*\omega_{\text{can}} +\sum_{1\leq i<j\leq
r}\Big((\text{pr}_i^*\omega_{\text{can}}
+\text{pr}_j^*\omega_{\text{can}})\cr &-{{\sqrt{-1}}\over
2}\sum_{k=1}^q(\text{pr}_i^*
\phi_k\wedge\text{pr}_j^*\bar\phi_k+\text{pr}_j^*\phi_k\wedge
\text{pr}_i^*\bar\phi_k)\Big),\cr}$$ which is well-known to be the pull-back
of the first Chern form of the Theta bundle together with the standard metric
induced by using theta norms. (See e.g., (4.3) below or [La2, p. 146].) Thus,
following the standard discussion in  [La2, Chapter VI, \S2-\S3], we complete
the proof of the theorem.
\vskip 0.30cm
\noindent
{\it Remark 4.1.1} It is not surprising that  the above 
theorem holds for a general $\omega$ as by definition we have the following
isometries:
$${\underline{{\Cal O}_M(Q_i)}}_\omega\|_{Q_j}
\simeq {\underline{{\Cal O}_M(Q_i)}}_{\omega_{\text{can}}}\|_{Q_j}(\simeq
{\underline{{\Cal O}_M(Q_i)}}_{\omega_{\text{can}}}|_{Q_j})
$$ and $${\underline{{\Cal O}_M(E)}}_\omega\|_{Q_i}
\simeq {\underline{{\Cal O}_M(E)}}_{\omega_{\text{can}}}\|_{Q_i}(\simeq
{\underline{{\Cal O}_M(E)}}_{\omega_{\text{can}}}|_{Q_i}).
$$ Similarly, we have 
$${\underline{K_M}}_\omega\|_{Q}
\simeq {\underline{K_M}}_{\omega_{\text{can}}}\|_{Q}(\simeq
{\underline{K_M}}_{\omega_{\text{can}}}|_{Q})$$ for any point $Q\in M$.
\vskip 0.30cm
\noindent
(4.2) In this section, we give further properties for the $\omega$-Faltings metrics.
\vskip 0.30cm
First of all, by definition, we have the following.
\vskip 0.30cm
\noindent
{\bf Lemma 4.2.1.} {\it With the same notation as above, there exists a
natural isometry
$$\Big(\lambda(K_M),h_{F;\omega}(\underline{K_M)}_\omega\Big)\simeq
\Big(\lambda(K_M),h_{F;\omega_{\text can}}(\underline{K_M}_{\omega_{\text
can}})\Big).\eqno(4.2.1)$$}
\vskip 0.20cm
On the other hand, by the proof of Theorem 4.1.1, we see that,  for general points 
$(Q_1,\dots,Q_q,Q)\in M^{q+1}$ 
such that $H^0(M,{\Cal O}_M(Q_1+\dots+Q_q-Q))=H^1(M,{\Cal
O}_M(Q_1+\dots+Q_q-Q))=\{0\}$,  
$\lambda({\Cal O}_M(Q_1+\dots+Q_q-Q))$ is simply ${\Bbb C}$, and the norm $1$
in ${\Bbb C}$ is proportional to
$\|\theta(Q_1+\dots+Q_r-Q)\|$, so that the ratio is independent of $(Q_1,\dots,Q_q,Q)$. Such a ratio gives an invariant
 associated to $(M,\omega)$. Following Faltings, we  define
the $\omega$-{\it Faltings delta function} $\delta(M,\omega)$ by
$$\|1\|_{h_{F;\omega}(\underline{{\Cal
O}_M(Q_1+\dots+Q_q-Q)})}=e^{-\delta(M;\omega)/8}\|
\theta(Q_1+\dots+Q_q-Q)\|.\eqno(4.2.2)$$
\vskip 0.30cm
\noindent
{\bf Proposition 4.2.2.} {\it With the same notation as above, we have 
$$\delta(M;\omega)=\delta(M;\omega_{\text{can}})(=\delta(M)).\eqno(4.2.3)$$
That is, $\omega$-Faltings delta function $\delta(M;\omega)$ is the same as
the original Faltings delta function $\delta(M)$.}
\vskip 0.30cm
\noindent
{\it Proof.} First of all, by (F3), for any point $Q\in M$, we have the
natural isometry
$$\Big(\lambda(K_M),h_{F;\omega}(\underline{K_M})\Big)\simeq 
\Big(\lambda(K_M(Q)),h_{F;\omega}(\underline{K_M(Q)})\Big)$$ due to the fact
that 
$$\underline{K_M(Q)}\|Q\simeq \underline{\Bbb C}.$$
Secondly, by a tedious calculation using (F3) again, we  get the following
isometry
$$\eqalign{~&\Big(\lambda(K_M(Q)),h_{F;\omega}(\underline{K_M(Q)})\Big)\cr
\simeq&
\Big(\lambda(K_M(Q-Q_1-\dots-Q_q)),
h_{F;\omega}(\underline{K_M(Q-Q_1-\dots-Q_q)})\Big)\cr
&\otimes\otimes_{i=1}^q\underline{K_M}\|_{Q_i}\otimes\otimes_{i=1}^q
\underline{{\Cal O}_M(Q)}\|_{Q_i}\cr &\otimes\otimes_{1\leq i<j\leq q}
(\underline{{\Cal O}_M(Q_i)}\|_{Q_j})^{\otimes -1}.\cr}$$
Thirdly, by (F4), the Serre isometry, we see that the last combination is
isometric to
$$\eqalign{~&\Big(\lambda({\Cal O}_M(Q_1+\dots+Q_q-Q)),
h_{F;\omega}(\underline{{\Cal O}_M(Q_1+\dots+Q_q-Q)})\Big)\cr
&\otimes\otimes_{i=1}^q\underline{K_M}\|_{Q_i}\otimes\otimes_{i=1}^q
\underline{{\Cal O}_M(Q)}\|_{Q_i}\cr
&\otimes\otimes_{1\leq i<j\leq q}
(\underline{{\Cal O}_M(Q_i)}\|_{Q_j})^{\otimes -1}.\cr}$$
Thus, we get the isometry
$$\eqalign{~&\Big(\lambda(K_M),h_{F;\omega}(\underline{K_M})\Big)\cr \simeq&
\Big(\lambda({\Cal O}_M(Q_1+\dots+Q_q-Q)),h_{F;\omega}(\underline{{\Cal
O}_M(Q_1+\dots+Q_q-Q)})\Big)\cr
&\otimes\otimes_{i=1}^q\underline{K_M}\|_{Q_i}\otimes\otimes_{i=1}^q
\underline{{\Cal O}_M(Q)}\|_{Q_i}\cr &\otimes\otimes_{1\leq i<j\leq q}
(\underline{{\Cal O}_M(Q_i)}\|_{Q_j})^{\otimes -1}.\cr}$$
Hence, finally, by Remark 4.1.1 and Lemma 4.2.1, we have the isometry
$$\eqalign{~&\Big(\lambda({\Cal
O}_M(Q_1+\dots+Q_q-Q)),h_{F;\omega}(\underline{{\Cal
O}_M(Q_1+\dots+Q_q-Q)}_\omega)\Big)\cr
\simeq&
\Big(\lambda({\Cal
O}_M(Q_1+\dots+Q_q-Q)),h_{F;\omega_{\text{can}}}(\underline{{\Cal
O}_M(Q_1+\dots+Q_q-Q)}_{\omega_{\text{can}}})\Big)\cr}$$ which, by
definition, completes the proof of the lemma.
\vskip 0.30cm
\noindent
{\it Remark 4.2.1.} We sometimes call Lemma 4.2.1 and Proposition 4.2.2 Mean
Value Lemmas too. 
\vskip 0.30cm
\noindent
(4.3) With the above definition of $\omega$-Faltings metric, we also have the
Noether isometry  without any further difficulty. Following Faltings [Fa] and
Moret-Bailly [MB], with arithmetic applications in mind, this is done as
follows. (I include this subsection simply for completeness which in turn offers me
a chance to give one of the main result of this paper, the $\omega$-Noether
isometry. If the reader does not want to waste his time on the known discussion
about theta norms, he  may simply jump to Theorem 4.3.1.)
\vskip 0.30cm
On the stack ${\Cal U}_q$ of principally polarized abelian varieties of dimension $q$, define a degree
$2^{2q}$ covering ${\Cal Q}$ which classifies pairs consisting of an abelian variety together with a symmetric ample
divisor which defines a principal polarization. Similarly, on the stack ${\Cal M}_q$ of regular 
algebraic curves of genus $q$, define the covering ${\Cal P}$, which
classifies pairs consisting of  regular genus $q$ curves together with one
of its theta-characters. Then over ${\Bbb Z}[1/2]$, by using the Abel-Jacobi
map, there  exists a cartesian diagram
$$\matrix {\Cal P}&\to&{\Cal Q}\cr
\downarrow&&\downarrow\cr
{\Cal M}_q&\to&{\Cal U}_q.\cr
\endmatrix$$
Over ${\Cal Q}$, there is the universal abelian variety 
$p:{\Cal A}\to {\Cal Q}$ together with the theta divisor 
$\Theta\subset{\Cal A}$, flat over ${\Cal Q}$. Denote the zero section of ${\Cal A}\to {\Cal Q}$ by $s$, 
then (up to a universal constant, over the corresponding analytic space,) we
have a natural  isometry $$s^*(\underline{\Omega_p^q})\simeq s^*
\underline{{\Cal O}_{\Cal A}(\Theta)}^{\otimes 2},$$ defined by multiplying $\theta^{-2}$.
Here $\underline{\Omega_p^q}$ denotes the line bundle $\Omega_p^q$ together
with the metric induced from the natural pairing  
$(\phi,\psi)\mapsto(-1)^{{{q(q-1)}\over 2}}
\cdot({{\sqrt{-1}}\over 2})^q\int_A\phi\wedge\bar\psi,$ while 
$\underline{{\Cal O}_{\Cal A}(\Theta)}$ denotes the line bundle ${\Cal O}_{\Cal A}(\Theta)$ together with the 
Hermitian metric defined by using the theta norm, i.e.,
$$\|\theta(Z,z)\|^2:=\sqrt{\text{det}Y}\cdot\exp(-2\pi^tyY^{-1}y)\cdot|\theta(Z,z)|^2$$ with $$\theta(Z,z):=
\sum_{n\in{\Bbb Z}^q}e^{\pi\sqrt{-1}^tnZn+2\pi\sqrt{-1}^tn\cdot z}.$$
Here a principally polarized abelian variety is taking of the form 
${\Bbb C}^q/({\Bbb Z}^q+Z\cdot{\Bbb Z}^q)$
for some complex $q\times q$ matrix $Z$ with positive definite imaginary part
$Y$.
\vskip 0.30cm
Thus if we denote the universal curve over ${\Cal P}$ by $p:{\Cal X}\to {\Cal P}$, 
there is a universal theta-character ${\Cal L}$ on ${\Cal X}$. In particular,
the $\omega$-Faltings metric on $\lambda({\Cal L})$ gives a Hermitian
metric on the associated line bundle over the analytic space corresponding
to
${\Cal P}$. On the other hand, we know that $\lambda({\Cal L})$ is simply the
pull-back of
$s^*{\Cal O}(-\Theta)$. So by definition, we see that such an isomorphism
gives an isometry for the corresponding Hermitian line  bundles if we
multiplies the Hermitian norm by
$\exp(\delta(M;\omega)/8)$  up to a constant depending only on the genus $q$
on each  connected component of ${\Cal P}$. 
\vskip 0.30cm
Finally, by (F0), we see that the Hermitian line bundle $\lambda(K_p)$ 
together with $\omega$-Faltings 
metric is merely the pull-back of the Hermitian line bundle $s^*\underline{\Omega_p^q}$. Thus, by using
Proposition 4.2.2, together with the detailed discussion of the constants
appeared above by Moret-Bailly in [MB],  we arrive at the following
\vskip 0.30cm
\noindent
{\bf Theorem 4.3.1.} ($\omega$-Noether isometry) {\it With respect to the
normalized volume $\omega$ (associated to a quasi-hyperbolic metric)  
on a compact Riemann surface $M$, for any $\omega$-admissible Hermitian line
bundle $\bar L$, we have the following isometry:
$$\Big(\lambda(L),h_{F;\omega}(\bar L)\Big)^{\otimes 12}\simeq\langle\bar
L,\bar L\otimes\underline{K_M}^{\otimes -1}\rangle^{\otimes 6}
\otimes\langle\underline{K_M},\underline{K_M}\rangle\otimes
{\Cal O}(e^{\delta(M)}\cdot(2\pi)^{-4q}).\eqno(4.3.1)$$}
\vskip 0.20cm
\noindent
{\it Remark 4.3.1.} The discussion in this section
works is simply due to the fact that  we here  use only the  normalized
$\omega$-Arakelov metrics on $K_M$ and ${\Cal O}_M(P)$.
\vskip 0.20cm
\noindent
{\it Remark 4.3.2.} Note that the Riemann-Roch isometry
gives the difference between  Hermitian line bundles 
$\Big(\lambda(K_M^{\otimes n}),h_{F;\omega}(\underline{K_M}^{\otimes
n})\Big)$ and 
$\Big(\lambda(K_M),h_{F;\omega}(\underline{K_M})\Big)$ for all half integers
$n$. So there are still some freedom for us to choose the
$\omega$-Riemann-Roch metric.  But the above discussion gives summation of
$\Big(\lambda(K_M^{\otimes 1/2}), h_{F;\omega}(\underline{K_M}^{\otimes
1/2})\Big)$ and 
$\Big(\lambda(K_M),h_{F;\omega}(\underline{K_M})\Big)$ for $\omega$-Faltings
metrics.  So we get a unique metric $h_{F;\omega}(\bar L)$ for all
$\omega$-admissible  Hermitian line bundles $\bar L$, so as to obtaining the
Noether isometry.  Such an idea was also previously used by Beilinson and
Manin in [BM] when they gave  Mumford  volume forms on  moduli spaces of
compact Riemann surfaces. It is for this reason we call (3.2.1) the weak
Riemann-Roch theorem, while we call the Deligne-Riemann-Roch theorem the
strong Riemann-Roch theorem. (See e.g., [We1, A1].)
\vskip 0.45cm
\noindent
{\we \S5. New metrics on determinant of cohomology
for singular metrics and Mumford type isometries}
\vskip 0.30cm
\noindent
(5.1)  We start with the following
\vskip 0.30cm
\noindent
{\bf Theorem 5.1.1.} (Deligne-Riemann-Roch Isometry for Singular Metrics) {\it 
For any normalized volume form $\omega$ on a compact Riemann surface $M$ associated
to a smooth metric or a quasi-hyperbolic metric, for any
$\omega$-admissible metric $\bar L$ on $M$, there exists an
$\omega$-determinant metric $h_{\overline{K_M}}(\bar L)$ on $\lambda(L)$
such that we have the following canonical
isometry
$$\Big(\lambda(L),h_{\overline{K_M}}(\bar L)\Big)^{\otimes
12}\simeq\langle\bar L,\bar L\otimes\overline{K_M}^{\otimes
-1}\rangle^{\otimes 6}
\otimes\langle\overline{K_M},\overline{K_M}\rangle\otimes
{\Cal O}(e^{a(q)}).\eqno(5.1.1)$$ Here $a(q):=(1-q)(24\zeta_{\Bbb Q}'(-1)-1)$
denotes the Deligne constant.}
\vskip 0.30cm
\noindent
{\it Proof.} First, let us prove the theorem when the metric on $K_M$
is simply the $\omega$-Arakelov metric, i.e., we for the time being assume that
$\overline{K_M}=\underline{K_M}$. In
this case,  by Theorem 4.3.1,  for any
$\omega$-admissible metric $\bar L$ on $M$, there exists a metric $h_{F,\omega}(\bar
L)$ on $\lambda(L)$ such that we have the following $\omega$-Noether isometry:
$$\Big(\lambda(L),h_{F;\omega}(\bar L)\Big)^{\otimes
12}\simeq\langle\bar L,\bar L\otimes\underline{K_M}^{\otimes
-1}\rangle^{\otimes 6}
\otimes\langle\underline{K_M},\underline{K_M}\rangle\otimes
{\Cal O}(e^{\delta(M)}\cdot(2\pi)^{-4q}).\eqno(5.1.2)$$
As a direct consequence,
if we set $$h_{\underline{K_M}}(\overline{L}):=
h_{F,\omega}(\overline{L})\cdot
e^{-\delta(M,\omega)/12}\cdot(2\pi)^{4q/12}\cdot e^{a(q)/12},\eqno(5.1.3)$$ with
$\delta(\omega,M)=\delta(M)$ the Faltings $\delta$-function, we then arrive at the
isometry
$$\Big(\lambda(L),h_{\underline{K_M}}(\bar L)\Big)^{\otimes 12}\simeq\langle\bar
L,\bar L\otimes\underline{K_M}^{\otimes -1}\rangle^{\otimes 6}
\otimes\langle\underline{K_M},\underline{K_M}\rangle\otimes
{\Cal O}(e^{a(q)}).\eqno(5.1.4)$$
\vskip 0.30cm
In general, we have $\overline {K_M}=\underline{K_M}\cdot {\Cal O}_M(e^{c})$ for a
certain constant function  $c$ on $M$ by applying the $\omega$-admissible
condition. So up to a constant $A(c,q,d)$, which can be easily
evaluated by using $\omega$-intersection and depends only on $c$
and $q$ and the degree $d$ of $L$,
 $$\langle\bar
L,\bar L\otimes\underline{K_M}^{\otimes -1}\rangle^{\otimes 6}
\otimes\langle\underline{K_M},\underline{K_M}\rangle$$
is simply $$\langle\bar L,\bar L\otimes\overline{K_M}^{\otimes
-1}\rangle^{\otimes 6}
\otimes\langle\overline{K_M},\overline{K_M}\rangle.$$ (We leave
the precise valuation of $A(c,q,d)$ to the reader as it is an interesting exercise
to understand the Polyakov variation formula for our metric.)
Set $$h_{\overline{K_M}}(\overline L):=h_{\underline{K_M}}(\overline L)\cdot
A(c,q,d),$$ 
 then we have the Deligne-Riemann-Roch isometry stated in the theorem. This
completes the proof of the theorem.
\vskip 0.30cm
As  direct consequences, we have the following.
\vskip 0.30cm
\noindent
{\bf Corollary 5.1.2.} (Mean Value Lemma) {\it With the same notation as above,
suppose that
$\omega_1$ and $\omega_2$ are two normalized volume forms on $M$, then
there exists a canonical isometry
$$\Big(\lambda(K_M), h_{\underline{K_M}_{\omega_1}}(\underline{K_M}_{\omega_1})\Big)
\simeq \Big(\lambda(K_M),
h_{\underline{K_M}_{\omega_2}}(\underline{K_M}_{\omega_2})\Big).$$}
\vskip 0.30cm
\noindent
{\it Proof.} By applying Theorem 5.1.1, the isometry is obtained by the Mean Value
Lemma in $\omega$-intersection theory, i.e., Proposition 2.5.3.
\vskip 0.30cm
For completeness, we give the following better than nothing
\vskip 0.30cm
\noindent
{\bf Corollary 5.1.3.} {\it With the same notation as above, assume that $\omega$ is
smooth on $M$, then $h_{\overline{K_M}}(\overline L)$ is the Quillen metric
$h_Q(\overline {K_M},\bar L)$ on $\lambda(L)$ associated to the metric on $M$
induced from $\overline {K_M}$ and the metrized line bundle $\bar L$.}
\vskip 0.30cm
\noindent
{\it Proof.} By applying Deligne-Riemann-Roch theorem  for
Quillen metric $h_Q(\overline{K_M},\bar L)$, (see e.g. [De2] together with [So],) we
have the isometry
$$\Big(\lambda(L),h_Q(\overline {K_M},\bar L)\Big)^{\otimes
12}\simeq\langle\bar L,\bar L\otimes\overline{K_M}^{\otimes
-1}\rangle^{\otimes 6}
\otimes\langle\overline{K_M},\overline{K_M}\rangle\otimes
{\Cal O}(e^{a(q)}).\eqno(5.1.2)$$ 
Comparing this with (5.1.4), we complete
the proof of this corollary.
\vskip 0.30cm
Moreover, from the above, easily, one sees that the determinant metric
$h_{\overline{K_M}}(\bar L)$ introduced in Theorem 5.1.1 is compactible with the
normalization process given in \S3. That is to say, we have the Polyakov
variation formula, the Mean Value Lemmas, among others.  Therefore, all
the above discussion is compactible. Surely,
{\it $h_{\overline{K_M}}(\bar L)$ is then the new metric on $\lambda(L)$ associated
to possibly singular metrics on ${\overline{K_M}}$ and on  $\bar L$ we seek at
the very beginning}. In the sequel, we give some applications and a more 
geometric interpretation of this new metric.
\vskip 0.30cm
\noindent
(5.2)  We start with a more suitable version of
Deligne-Riemann-Roch isomorphism for punctured Riemann surfaces. 
\vskip 0.30cm
First,  the algebraic Noether theorem tells us that for a compact Riemann
surface $M$ there is a canonical isomorphism $$\lambda({\Cal O}_M)^{\otimes
12}\simeq \langle K_M,K_M\rangle.$$ Secondly, the adjunction isomorphism gives the
following isomorphisms
$$\langle K_M(P_i),P_i\rangle\simeq {\Cal O}.$$ 
Here $P_i$, $i=1,\dots,N,$ denotes the punctures of $M$. As a direct consequence, we
have
$$\eqalign{~&\langle K_M,K_M\rangle\cr
\simeq&\langle
K_M(P_1+\dots+P_N),K_M(P_1+\dots+P_N)\rangle\cr
&\otimes\Big(\langle K_M,{\Cal
O}_M(P_1+\dots+P_N)\rangle\otimes
\langle {\Cal O}_M(P_1+\dots+P_N), K_M(P_1+\dots+P_N)\rangle\Big)^{-1}\cr
\simeq&
\langle K_M(P_1+\dots+P_N),K_M(P_1+\dots+P_N)\rangle\cr
&\otimes\Big(\langle K_M,{\Cal
O}_M(P_1+\dots+P_N)\rangle\otimes \otimes_{i=1}^N\langle{\Cal
O}_M(P_i),K_M(P_i)\otimes \otimes_{j=1,\not=i}^N{\Cal O}_M(P_j)\rangle\Big)^{-1}\cr
\simeq&\langle
K_M(P_1+\dots+P_N),K_M(P_1+\dots+P_N)\rangle\cr
&\otimes\Big(\langle K_M,{\Cal
O}_M(P_1+\dots+P_N)\rangle\otimes \otimes_{i=1}^N\langle{\Cal
O}_M(P_i), \otimes_{j=1,\not=i}^N{\Cal O}_M(P_j)\rangle\Big)^{-1}\cr
\simeq&
\langle K_M(P_1+\dots+P_N),K_M(P_1+\dots+P_N)\rangle\cr
&\otimes\Big(\langle K_M,{\Cal
O}_M(P_1+\dots+P_N)\rangle\otimes \otimes_{1\leq i<j\leq N}\langle{\Cal
O}_M(P_i), {\Cal O}_M(P_j)\rangle^{\otimes 2}\Big)^{-1}.\cr}$$ Thus if we set
$$\Delta_\alpha=\cases \langle
K_M(P_1+\dots+P_N),K_M(P_1+\dots+P_N)\rangle,&$if\ $\alpha=0$$,\\
\otimes_{k=1}^N\langle {\Cal
O}_M(P_k),{\Cal
O}_M(P_k)\rangle(=\otimes_{k=1}^N\langle K_M,{\Cal
O}_M(P_k)\rangle^{\otimes -1}),&$if\ $\alpha=1$$,\\
\otimes_{1\leq i<j\leq N}\langle{\Cal
O}_M(P_i), {\Cal O}_M(P_j)\rangle,&$if\ $\alpha=2$$.\endcases\eqno(5.2.1)$$ 
we have the following canonical isomorphism $$\lambda({\Cal O}_M)^{\otimes
12}\simeq \Delta_0\otimes
\Delta_1\otimes\Delta_2^{\otimes -2}.\eqno(5.2.2)$$ 
In this way, we arrive at the following version
of Deligne-Riemann-Roch isomorphism, which is the most suitable one for punctured
Riemann surfaces.
\vskip 0.30cm
\noindent
{\bf Theorem 5.2.1.} (Deligne-Riemann-Roch Theorem) {\it With the
same notation as above, for all line bundles $L$ on $M$, we have
$$\lambda(L)^{\otimes 12}\simeq
\Delta_0\otimes\Delta_1\otimes\Delta_2^{\otimes -2}\otimes\langle L,L\otimes
K_M^{\otimes -1}\rangle^{\otimes 6}.\eqno(5.2.3)$$}
\vskip 0.30cm
\noindent
{\it Remark 5.2.1.} The reader may wonder why we use
$\Delta_0$ as a very basic object to build up the
isomorphism. We here justify our choice by the following two reasons: first
of all, the logarithmic geometry says that for punctured Riemann surfaces,
$K_M(P_1+\dots+P_N)$ is far more natural (see e.g., [Fu]); second, we
will see that on the Teichm\"uller space of the punctured Riemann surfaces of
signature $(q,N)$, the most natural line bundle corresponding to the Weil-Petersson
K\"ahler form is given by the Deligne pairing $\langle
K_M(P_1+\dots+P_N),K_M(P_1+\dots+P_N)\rangle$, rather than
$\langle K_M,K_M\rangle$.
\vskip 0.30cm
To go further, for punctures Riemann surfaces $M^0$ with cusps $P_1,\dots,P_N$ and
the smooth compactification $M$,  we define the Mumford type line bundles
$\lambda_n$ by
$$\lambda_n:=\cases\lambda(K_M^{\otimes n}\otimes({\Cal
O}_M(P_1+\dots+P_N))^{\otimes n-1}),&$\quad\text{if}\ $n>0$$;\\
\lambda({\Cal O}_M),&$\quad\text{if}\ $n=0$$;\\
\lambda((K_M(P_1+\dots+P_N))^{\otimes
n}),&$\quad\text{if}\ $n<0$$.\endcases\eqno(5.2.4)$$
\vskip 0.20cm
\noindent
{\it Remark 5.2.2.} For the time being, we justify
 this definition of Mumford type line bundles $\lambda_n$ as follows.
First of all, the most natural line bundle associated to a punctured Riemann
surface is the associated logarithmic tangen line bundle, so it is fairly natural
to define $\lambda_n$ for $n$ negative by setting 
$$\lambda_n:=\lambda((K_M(P_1+\dots+P_N))^{\otimes n}).$$ Second, Serre duality
should give intrinsic relations among all $\lambda_n$'s. This then gives the above
defintion of $\lambda_n$ for $n$ positive. 
\vskip 0.30cm
\noindent
{\bf Theorem 5.2.2.} (Generalized Mumford Relations) {\it With the same
notation as above, for all positive integers $n$, we have the following
isomorphisms:
\vskip 0.30cm
\noindent
(a) $\lambda_n\simeq\lambda_{1-n};$ 
\vskip 0.30cm
\noindent
(b) $\lambda_n^{\otimes 12}\simeq \Delta_0^{\otimes
(6n^2-6n+1)}\otimes \Delta_1\otimes\Delta_2^{\otimes 10-12n},$ and
\vskip 0.30cm
\noindent
(c) $\lambda_n\simeq \lambda_0^{\otimes(6n^2-6n+1)}\otimes 
\Delta_1^{\otimes -{{n(n-1)}\over 2}}\otimes\Delta_2^{\otimes(n-1)^2}.$}
\vskip 0.30cm
\noindent
{\it Proof.} (a) is a direct consequence of the definition. The proofs of (b) and
(c) are similar, so we only give the one for (b).
For this latest purpose, we have the following calculation. In Theorem 5.2.1,
setting $L:=K_M^{\otimes n}\otimes({\Cal
O}_M(P_1+\dots+P_N))^{\otimes n-1}$, we get
$$\eqalign{~&\lambda_n^{\otimes 12}\cr
\simeq&\Delta_0\otimes\Delta_1\otimes
\Delta_2^{\otimes -2}\otimes\langle K_M^{\otimes n}\otimes({\Cal
O}_M(P_1+\dots+P_N))^{\otimes n-1},\cr
&\hskip 5.0cm K_M^{\otimes n}\otimes({\Cal
O}_M(P_1+\dots+P_N))^{\otimes n-1}\otimes
K_M^{\otimes -1}\rangle^{\otimes 6}\cr
\simeq&\Delta_0\otimes\Delta_1\otimes
\Delta_2^{\otimes -2}\cr
&\qquad \otimes\langle K_M^{\otimes n}\otimes({\Cal
O}_M(P_1+\dots+P_N))^{\otimes n-1},K_M(P_1+\dots+P_N)\rangle^{\otimes
6(n-1)}\cr
\simeq&\Delta_0^{\otimes (6n^2-6n+1)}\otimes\Delta_1\otimes
\Delta_2^{\otimes -2}\cr
&\qquad \otimes\langle
{\Cal O}_M(P_1+\dots+P_N),K_M(P_1+\dots+P_N)\rangle^{\otimes -6(n-1)}.\cr}$$ Hence,
by the adjunction isomorphism, as in the proof of (5.2.2), we have
$$\lambda_n^{\otimes 12}\simeq \Delta_0^{\otimes (6n^2-6n+1)}\otimes\Delta_1\otimes
\Delta_2^{\otimes -12n+10}.$$ This completes the proof of the theorem.
\vskip 0.30cm
\noindent
(5.3) Now we give the counter part of the metric theory for the discussion 
in (5.2). We start with a discussion on results in which only $\omega$-Arakelov
metrics on both the canonical line bundle and pointed line bundles associated to
cusps are used.
\vskip 0.30cm
For a normalized volume form $\omega$ on $M$, define the following metrized
lines:
$$\underline{\lambda_n}:=\cases \Big(\lambda_n,
h_{\underline{K_M}}(\underline{K_M}^{\otimes n}\otimes {\Cal O}_M(
\underline{P_1+\dots+P_N})^{\otimes n-1})\Big),&$\ \text{if}\ $n>0$$;\\
\Big(\lambda_0,h_{\underline{K_M}}(\underline{{\Cal O}_M}),&$\ \text{if}\
$n=0$$;\\
\Big(\lambda_n,
h_{\underline{K_M}}((\underline{K_M}(
\underline{P_1+\dots+P_N}))^{\otimes n}\Big),&$\ \text{if}\
$n<0$$.\endcases\eqno(5.3.1)$$
$$\underline{\Delta_n}:=\cases\langle
\underline{K_M}(\underline{P_1+\dots+P_N}),\underline{K_M}(
\underline{P_1+\dots+P_N})\rangle,&$\ \text{if}\ $n=0$$;\\
\otimes_{k=1}^N\langle
\underline{{\Cal O}_M(P_k)},\underline{{\Cal O}_M(P_k)}\rangle,&$\ \text{if}\
$n=1$$;\\
\otimes_{1\leq i<j\leq N}\langle
\underline{{\Cal O}_M(P_i)},\underline{{\Cal O}_M(P_j)}\rangle,&$\ \text{if}\
$n=2$$.\endcases\eqno(5.3.2)$$
 Note that in this case the adjunction isometry holds, so the proof for the
algebraic results in (5.2) is valid here without any essential change. In other
words, we have the  following
\vskip 0.30cm
\noindent
{\bf Theorem 5.3.1.} {\it With the same notation as above, for any positive
integer $n$, we have the following isometries:
\vskip 0.30cm
\noindent
(a) (Serre isometry) 
$$\underline{\lambda_n}\simeq\underline{\lambda_{1-n}};$$
\vskip 0.30cm
\noindent
(b) (Generalized Mumford isometry) 
$$\underline{\lambda_n}^{\otimes 12}\simeq
\underline{\Delta_0}^{\otimes 6n^2-6n+1}\otimes
\underline{\Delta_1}\otimes\underline{\Delta_2}^{\otimes
-12+10}\otimes{\Cal O}(e^{a(q)});$$
\vskip 0.30cm
\noindent
(c) (Generalized Mumford isometry) 
$$\underline{\lambda_n}\simeq
\underline{\lambda_1}^{\otimes 6n^2-6n+1}\otimes
\underline{\Delta_1}^{\otimes -{{n(n-1)}\over 2}}
\otimes\underline{\Delta_2}^{\otimes(n-1)^2}\otimes{\Cal O}
(e^{-{{n(n-1)}\over 2}\cdot a(q)}).$$}
\vskip 0.20cm
\noindent
{\it Proof.} As said before, the proof of this theorem is
essentially given in (5.2), as we here use the
$\omega$-Arakelov metrics so that 
the adjunction isometry holds. With this, by using Theorem 5.1.1 with a tedious
calculation, we complete the proof of this theorem.
\vskip 0.30cm
\noindent
(5.4) More generally, without using the adjunction isometry and with the application
to the moduli problems in mind, we in this subsection give a generalization of
Theorem 5.3.1. As in (5.3), we always fix a normalized volume form $\omega$ on
$M$.
\vskip 0.30cm
For an $n+1$-tuple of real numbers $(\alpha;\beta_1,\dots,\beta_N)$,
define the associated metrized lines as follows:
$$\overline{\lambda_n}^{\alpha;\beta}:=\cases\Big(\lambda_n,
h_{\overline{K_M}^\alpha}((\overline{K_M}^\alpha)^{\otimes n}\otimes {\Cal
O}_M(
\overline{P_1}^{\beta_1}+\dots+\overline{P_N}^{\beta_N})^{\otimes
n-1})\Big),&$\ \text{if}\ $n>0$$;\\ 
\Big(\lambda_0,h_{\overline{K_M}^\alpha}(\underline{{\Cal
O}_M})\Big),&$\ \text{if}\ $n=0$$;\\
\Big(\lambda_n,
h_{\overline{K_M}^\alpha}((\overline{K_M}^\alpha(
\overline{P_1}^{\beta_1}+\dots+\overline{P_N}^{\beta_N}))^{\otimes n})\Big),&$\
\text{if}\ $n<0$$;\endcases\eqno(5.4.1)$$ and
$$\overline{\Delta_n}^{\alpha;\beta}:=\cases\langle
\overline{K_M}^\alpha(
\overline{P_1}^{\beta_1}+\dots+\overline{P_N}^{\beta_N}),
\overline{K_M}^\alpha(
\overline{P_1}^{\beta_1}+\dots+\overline{P_N}^{\beta_N})
\rangle,&$\ \text{if}\ $n=0$$;\\
\langle \overline{K_M}^\alpha,{\Cal O}_M
(\overline{P_1}^{\beta_1}+\dots+\overline{P_N}^{\beta_N})
\rangle^{\otimes -1},&$\ \text{if}\
$n=1$$;\\
\langle\overline{K_M}^\alpha(
\overline{P_1}^{\beta_1}+\dots+\overline{P_N}^{\beta_N}),
{\Cal O}_M
(\overline{P_1}^{\beta_1}+\dots+\overline{P_N}^{\beta_N})\rangle^{\otimes
{1\over 2}},&$\
\text{if}\ $n=2$$.\endcases\eqno(5.4.2)$$ Here
$\overline{K_M}^{\alpha}:=\underline{K_M}\otimes
{\Cal O}_M(e^\alpha)$ and
$\overline{P_i}^{\beta_i}:=\underline{P_i}\otimes {\Cal O}_M(e^{\beta_i}),$
$i=1,\dots,N$. Then we get the following
\vskip 0.30cm
\noindent
{\bf Theorem 5.4.1.} {\it With the same notation as above, for any positive
integer $n$, we have the following isometries:
\vskip 0.30cm
\noindent
(a) (Serre isometry) 
$$\overline{\lambda_n}^{\alpha;\beta}\simeq\overline{\lambda_{1-n}}^{\alpha;\beta};$$
\vskip 0.30cm
\noindent
(b) (Generalized Mumford isometry) 
$$(\overline{\lambda_n}^{\alpha;\beta})^{\otimes 12}\simeq
(\overline{\Delta_0}^{\alpha;\beta})^{\otimes 6n^2-6n+1}\otimes
(\overline{\Delta_1}^{\alpha;\beta})
\otimes(\overline{\Delta_2}^{\alpha;\beta})^{\otimes
-12+10}\otimes{\Cal O}(e^{a(q)});$$
\vskip 0.30cm
\noindent
(c) (Generalized Mumford isometry) 
$$\overline{\lambda_n}^{\alpha;\beta}\simeq
(\overline{\lambda_1}^{\alpha;\beta})^{\otimes 6n^2-6n+1}\otimes
(\overline{\Delta_1}^{\alpha;\beta})^{\otimes -{{n(n-1)}\over 2}}
\otimes(\overline{\Delta_2}^{\alpha;\beta})^{\otimes(n-1)^2}\otimes{\Cal
O}(e^{-{{n(n-1)}\over 2}\cdot a(q)}).$$}
\vskip 0.20cm
\noindent
{\it Proof.} By Proposition 3.2.3, i.e., the Serre duality, we have (a). The proofs
of (b) and (c) come directly from applying the Deligne-Riemann-Roch isometry for
singular metrics, i.e., Theorem 5.1.1. Indeed, we,  in the above definition (5.4.2)
on
$\overline{\Delta_1}^{\alpha;\beta}$ and $\overline{\Delta_2}^{\alpha;\beta}$, 
already made a subtle  change from these in (5.3.2). Moreover, it is clear  by a
direct calculation that with (5.4.2) no adjunction isometry is needed to get (b)
and (c). We leave the details to the reader.
\vskip 0.30cm
\noindent
{\we Appendix to \S5. Universal Riemann-Roch Isomorphism}
\vskip 0.30cm
\noindent
Even though the above discussion is about a single punctured
Riemann surface,  but it is clear that it works for a family. For our own
convienence (to the later discussion) and  for the
completeness, we in this appendix to \S5 study briefly the structure of the moduli
space of punctured Riemann surfaces and their compactifications, which is
well-known to experts. The details for most of the
statements, except for the universal Riemann-Roch theorem for
punctured curves,  may be found in [DM] and Kundsen's series of papers ([KM]
and [Kn]).  
\vskip 0.30cm
Denote the moduli space of $N$ ordered  pointed compact Riemann
surfaces of genus $q$ by ${\Cal M}_{q,N}$, and its universal curve by
$\pi_{q,N}: {\Cal C}_{q,N}\to {\Cal M}_{q,N}$. (Here we for
simplicity assume that
${\Cal M}_{q,N}$ is a fine moduli. If the reader does not like this,
she or he should then use  $V$-manifold language or algebraic
stack language to justify what follows.)
\vskip 0.30cm
Naturally, we have the relative canonical line bundle $K_{q,N}$ for
$\pi_{q,N}$, and $N$-sections ${\Bbb P}_1,\dots,{\Bbb P}_N$
corresponding to $N$ points on each surface. One knows that ${\Cal
C}_{q,N}$ may be viewed as ${\Cal M}_{q,N+1}$, while $\pi_{q,N}$ is
simply the map of dropping the last, i.e., the $(N+1)$-th point. In
particular, we have the following commutative diagram:
$$\matrix {\Cal M}_{q,N+1}=&{\Cal C}_{q,N}&\buildrel
\phi_{q,N}\over\to&{\Cal C}_{q,N-1}&={\Cal M}_{q,N}\\
&&&&\\
&\pi_{q,N}\downarrow&&\downarrow\pi_{q,N-1}&\\
&&&&\\
&{\Cal M}_{q,N}&\buildrel\pi_{q,N-1}\over\to& {\Cal
M}_{q,N-1}&.\endmatrix$$ Here $\phi_{q,N}$ viewed as a morphism from
${\Cal M}_{q,N+1}$ to ${\Cal M}_{q,N}$ is simply the morphism defined
by dropping the second to the last point.
\vskip 0.30cm
To compactify ${\Cal M}_{q,N}$, we need to add two types of
boundaries. That is, the boundaries coming from the degeneration of
compact Riemann surfaces of genus $q$, and the boundaries coming from
the degeneration of punctures. For our own convinence, we call the
first type of boundaries  the absolute horizontal boundaries, while
we call the second type of boundaries  the relative horizontal
boundaries.
\vskip 0.30cm
For  absolute horizontal boundaries, it is well-known that as a
codimension one subvarity, it consists of
$[{q\over 2}]+1$ irreducible components
$\Delta_0,\Delta_1,\dots,\Delta_{[{q\over 2}]}$. (Please do confuse
$\Delta$ here with $\Delta$ elsewhere in this paper.) Indeed, such  boundaries may
be at best understood via the universal curve: for a general point
$x\in
\Delta_0$, the corresponding fiber is a genus
$q$ curve with one non-separating node; while for
a general point $x\in \Delta_i$, $i=1,\dots,[{q\over 2}]$, the
corresponding fiber is a genus
$q$ curve with one separating node whose only two irreducible
components are smooth and of genera $i$ and $q-i$ respectively.
\vskip 0.30cm
To understand the relative horizontal boundaries, we suggest the
reader to consult  papers on resolution of diagonals, say, [BG] and
[FM]. In any case, this may be understood by considering what happens
for a general compact Riemann surface. So we for the time being
assume that $M$ is a compact Riemann surface without non-trivial
automorphisms. Then the moduli $M^{(1)}$ of a point on $M$ is simply
$M$ itself. And the universal curve over $M^{(1)}$ is simply
$M\times M$ with the section given by the diagonal. $M\times M$ can
be viewed as the moduli $M^{(2)}$ of an ordered pairing on $M$. There
is no problem to see the fiber together with two sections of the
universal curve over  points in $M^{(2)}$ which are away from the
diagonal: the fiber is simply
$M$ itself together with two distinct points (as we assume that the
base point in $M^{(2)}$ is off diagonal). On the other hand, if the
point in $M^{(2)}$ is on the diagonal, two sections intersect each
other. So  we cannot simply find two distinct points on $M$, the fake
fiber in the universal curve. To remedy this,
Grothendieck-Mumford-Knudsen first blow up $M$ at this point so that
two points can be pulled apart. In other words, the fibers over a
point on the diagonal now admits two irreducible components: one is
the original curve, while the other is a projective line together
with three marked points -- the intersection piont with $M$
representing the center of the blowing up, while the other two points
representing two infinitesimal pionts over the intersection of $M$
with
${\Bbb P}^1$. In this way, in particular, we see that the universal
curve admits two sections which can never meet each other.
\vskip 0.30cm
This picture may be generalized to the moduli $M^{(N)}$ of ordered $N$
points for $M$. To describe it together with its universal curve, we
 consider the set $\{1,\dots,N\}$. For each subset $S$ of
$\{1,\dots,N\}$ with cardinal number $\#S$ at least two, we have an
$S$-diagonal $D_S$ in $M^{N}$. We know that to have $M^{(N)}$, we
need to blow up these diagonals, so that we then get normal crossing
divisors $\Delta_S$ resulting from these diagonals $D_S$ (and
exceptional divisors). In particular, for a general point $x\in
\Delta_S$, the fiber of the universal curve consists of two
irreducible components, one is the original curve $M$ while the other
is the projective line ${\Bbb P}^1$. Moreover, on $M$ there are
$N-\#S$ marks and on ${\Bbb P}^1$, there are remaining $\#S$-marks.
Similarly, as for the case when $N=2$, now on the universal curve, there
are $N$-sections which do not intersect pairwise.
\vskip 0.30cm
From the above discussion, we see that the absolute
horizontal boundaries consist of divisors
$\Delta_0,\Delta_1,\dots,\Delta_{[{q\over 2}]}$, and the relative
horizontal boundaries consist of divisors $\Delta_S$ for each subset
$S$ of $\{1,\dots,N\}$ with cardinal number at least two. 
\vskip 0.30cm
\noindent
{\bf Universal Riemann-Roch Theorem.} {\it With the
same notation as above, denote by $\bar\pi_{q,N}:\overline{{\Cal
C}_{q,N}}\to \overline{{\Cal M}_{q,N}}$ the universal curve over
stably compactified moduli space of punctured Riemann surfaces of
signature $(q,N)$, and by $K$ the associated relative
canonical line bundle, then for any line bundle $L$ on $\overline{{\Cal
C}_{q,N}}$, we have the canonical isomorphism
$$\lambda(L)^{\otimes 12}\simeq
\langle L,L\otimes K\rangle^{\otimes
6}\otimes \langle K,K\rangle\otimes\otimes_{i=0}^{[{q\over 2}]}{\Cal
O}(\Delta_i)\otimes\otimes_{S\subset\{1,\dots,N\},\#S\geq 2}{\Cal
O}(\Delta_S).$$}
\vskip 0.30cm
The proof may be given by using Mumford's arguement on 
Riemann-Roch theorem for the universal curve over stably moduli space
of compact Riemann surfaces, as only ordinary double points are
involved here. We leave this to the reader. (See however [We2].)
\vskip 0.30cm
To end this appendix to section 5, we give the following list of relations for line
bundles associated to $\bar\pi_{q,N}$ coming from the intersection.
\vskip 0.30cm
\noindent
(a) $\displaystyle{\langle {\Bbb P}_i,{\Bbb P}_j\rangle\simeq {\Cal
O}}$, if $i,j=1,\dots,N$ and $i\not=j$;
\vskip 0.30cm
\noindent
(b) $\displaystyle{\langle K({\Bbb P}_i),{\Bbb P}_i\rangle\simeq {\Cal
O}}$, if $i=1,\dots,N$;
\vskip 0.30cm
\noindent
(c) $\displaystyle{\langle K,{\Bbb P}_i\rangle\simeq
\Big(\pi_{q,N-1}^*\langle K,{\Bbb P}_i\rangle\Big)({\Bbb P}_i)}$, if
$i+1,\dots,N-1$;
\vskip 0.30cm
\noindent
(d) $\displaystyle{\langle K,{\Bbb P}_N\rangle\simeq K({\Bbb
P}_1+\dots+{\Bbb P}_{N-1})}$;
\vskip 0.30cm
\noindent
(e) $\displaystyle{K({\Bbb
P}_1+\dots+{\Bbb P}_N)=\phi^*_{q,N}\Big(
K({\Bbb
P}_1+\dots+{\Bbb P}_{N-1})\Big)}.$
\vskip 0.30cm
Indeed, with these relations, we may easily generate many interesting
relations such as the dilation equation and the string equation. We
leave all this to the reader.
\vskip 0.30cm
\noindent
{\we \S6. Arakelov-Poincar\'e volume and a geometric interpretation of our new
metrics}
\vskip 0.30cm
\noindent
(6.1) From now on, we will apply our admissible theory to singular hyperbolic
metrics. For doing so, we need to understand how the
geometrically defined hyperbolic metric on the logarithmic tangent line
bundle relates with the arithmetically defined hyperbolic-Arakelov metric on the
canonical line bundle. We bridge them via an invariant, the Arakelov-Poincar\'e
volume, for a punctured Riemann surface. We end this section with  a geometric
interpretation for our new metric on determinant of cohomology. 
\vskip 0.30cm
Let us start with a discussion on hyperbolic metrics on punctured Riemann
surfaces. As before, denote by $\omega_{\text{hyp}}$ the normalized volume form
associated to the  standard hyperbolic metric
$\tau^0_{\text{hyp}}$ on a punctured Riemann surface $M^0$ of signature
$(q,N)$. Thus, in particular, if we denote the corresponding volume form
(with respect to $\tau^0_{\text{hyp}}$) by $d\mu_{\text{hyp}}$, then
$\int_{M^0}d\mu_{\text{hyp}}=2\pi(2q-2+N)$, and
$2\pi(2q-2+N)\omega_{\text{hyp}}=d\mu_{\text{hyp}}.$
\vskip 0.30cm
For $\tau^0_{\text{hyp}}$, or equivalently for $d\mu_{\text{hyp}}$ on $M^0$,
if we view them as a singular metric on $M$, the smooth compactification of $M^0$,
then the natural line bundle we should attach to it is the so-called
logarithmic tangent bundle
$T_M\langle\log D\rangle$. Here $D$ denotes the divisr at infinity, or the cuspidal
divisor, i.e., $P_1+\dots+P_N$. (See e.g., [De1], [Mu1] or [Fu]). Over the compact
Riemann surface $M$,
$T_M\langle\log D\rangle$ is nothing but the dual of the line bundle
$K_M(P_1+\dots+P_N)$. Here as before
$K_M$ denotes the canonical line bundle of $M$. So if we denote the induced 
Hermitian metric from $\tau^0_{\text{hyp}}$ on $K_M(P_1+\dots+P_N)$ by
$\tau_{\text{hyp};K_M(D)}^\vee$, we get the following Einstein equation
$$c_1\Big(K_M(P_1+\dots+P_N),\tau_{\text{hyp};K_M(D)}^\vee\Big)
=d\mu_{\text{hyp}}=(2q-2+N)\omega_{\text{hyp}}.\eqno(6.1.1)$$
\vskip 0.20cm
{\it We are not quite satisfied with this}, as the metric discussed above only has
its nice meaning on the logarithmic tangent bundle. In particular, it does not give
us any indication on how to get an $\omega_{\text{hyp}}$-admissible metric on $K_M$,
without which we cannot apply our admissible theory.
So we should seek new admissible metrics $\rho_{\text{hyp};K_M}$ and
$\rho_{\text{hyp};P_i}$ on $K_M$ and on  ${\Cal O}_M(P_i)$,
$i=1,\dots,N$, respectively, which naturally come  from the standard hyperbolic
metric. More precisely, for the time being, the picture we
have in mind for these admissible  metrics is that they are very natural in
the following sense:
\vskip 0.30cm
\noindent
(i) they should be $\omega_{\text{hyp}}$-admissible; 
\vskip 0.30cm
\noindent
(ii) they should give the following decomposition of the standard hyperbolic
metric   on $K_M(P_1+\dots+P_N)$
$$\tau_{\text{hyp};K_M(D)}^\vee=\rho_{\text{hyp};K_M}\otimes\rho_{\text{hyp};P_1}\otimes\dots\otimes
\rho_{\text{hyp};P_N};\eqno(6.1.2)$$
\vskip 0.30cm
\noindent
(iii) they should obey the residue isometry, i.e.,
$$(K_M(P_i),\rho_{\text{hyp};K_M}\otimes\rho_{\text{hyp};P_i})\|_{P_i}\simeq
\underline{\Bbb C}\eqno(6.1.3)$$ for all $i=1,\dots,N$.
\vskip 0.30cm
\noindent
{\it Remark 6.1.1.} There is  an interesting metric on $K_M$ induced from
$d\mu_{\text{hyp}}$, i.e., view it as the dual of the tangent bundle $T_M$,
then the singular volume form
$d\mu_{\text{hyp}}=g(z){{\sqrt{-1}}\over 2}dz\wedge d\bar z$ will give a 
singular metric via the function $g$ on $T_M$ and hence on $K_M$. But this
metric on $K_M$ is unfortunately not $\omega_{\text{hyp}}$-admissible. (Otherwise,
the problem should be much easier.) 
\vskip 0.30cm
Before defining the above metrics on $K_M$ and on ${\Cal O}_M(P_i)$,
$i=1,\dots N$, respectively, motivated by our previous work on admissible theory for
smooth hyperbolic metrics in [We1], we now introduce an invariant
$A_{\text{Ar,hyp}}(M^0)$, the Arakelov-Poincar\'e volume, associated to a
punctured Riemann surface $M^0$ as follows.
\vskip 0.30cm
First of all, following Selberg, define the so-called Selberg zeta function
$Z_{M^0}(s)$ of $M^0$ for $\text{Re}(s)>1$ by the absolutely convergent
product
$$Z_{M^0}(s):=\prod_{\{l\}}\prod_{m=0}^\infty(1-e^{-(s+m)|l|}),\eqno(
6.1.4)$$ where
$l$ runs over the set of all simple  closed geodesics on $M^0$ with
respect to the hyperbolic metric $d\mu_{\text {hyp}}$ on $M^0$, and $|l|$
denotes the length of $l$. It is known that by using  Selberg trace
formula for weight zero  forms the function $Z_{M^0}(s)$ admits a
meromorphic continuation to the whole complex
$s$-plane which has a simple zero at $s=1$. Secondly, motivated by the
work of D'Hoker-Phong ([D'HP]) and Sarnak ([Sa]), we introduce the
following factorization for the Selberg zeta function:
$$Z_{M^0}(s)=:\text{det}(\Delta_{\text {hyp}}+s(s-1))\cdot {\Bbb
N}(s)^{2q-2+N}.\eqno(6.1.5)$$ Here $\Delta_{\text {hyp}}$ denotes the
hyperbolic Laplacian on $M^0$, ${\Bbb N}(s)$ denotes the function
$${\Bbb
N}(s):={{e^{-E+s(s-1)}}\over{{2\pi}^s}}\cdot{{\Gamma(s)}\over{(\Gamma_2(s))^2}}
\eqno(6.1.6)$$
with $E=-{1\over 4}-{1\over 2}\log 2\pi+2\zeta_{\Bbb Q}'(-1)$, $\Gamma(s)$
the ordinary gamma function, and $\Gamma_2(s)$ the Barnes double gamma
funtion. Thirdly, define the regularized determinant for the Laplacian
$\Delta_{\text{hyp}}$ by
$$\text{det}^*(\Delta_{\text{hyp}}):=
{d\over{ds}}\Big(\text{det}(\Delta_{\text{hyp}}+s(s-1))\Big)\Big|_{s=1}.\eqno(6.1.7)$$
Finally,  define the {\it Arakelov-Poincer\'e volume}
$A_{\text{Ar,hyp}}(M^0)$ for
$M^0$ via the formula:
$$\log A_{\text{Ar,hyp}}(M^0):={a_{\text{Ar,hyp}}}:={{12}\over 2}\cdot
{1\over {2q-2}}\cdot\Big(\text{log}
{{\text{det}^*\Delta_{\text{Ar}}}\over {A_{\text{Ar}}(M)}}
-\text{log}{{\text{det}^*\Delta_{\text{hyp}}}\over
{2\pi(2q-2)}}\Big).\eqno(6.1.8)$$ 
Here $\Delta_{\text{Ar}}$ denotes the Laplacian for the
Arakelov metric on $M$, $A_{\text{Ar}}(M)$ denotes the volume of
$M$ with respect to the Arakelov metric. 
\vskip 0.30cm
\noindent
{\it Remark 6.1.2.} By definition, we know that, up to a constant factor depending
only on the signature $(q,N)$ of $M^0$, $\text{det}^*(\Delta_{\text{hyp}})$ is
simply $Z'_{M^0}(1)$. We leave this interesting point to the reader.
 Please also carefully compare our definition
(6.1.7) of the regularized determinant for the Laplacian with the one proposed by
Efrat in the one page correction of [Ef]. 
\vskip 0.30cm
\noindent
{\it Remark 6.1.3.} Obviously, the Arakelov-Poincar\'e volume is  a very
natural invariant for the punctured Riemann surface $M^0$, hence can be 
viewed as a certain interesting function on the Teichm\"uller space
$T_{q,N}$ of punctured Riemann surfaces of signature $(q,N)$. The reader may
consult [We1] for the degeneration behavior of this invariant when $N=0$. 
\vskip 0.30cm
\noindent
(6.2) With above discussion on the Arakelov-Poincar\'e volume for $M^0$, we
are ready to introduce the geometrically natural admissible metrics on $K_M$ and
${\Cal O}_M(P_i)$, $i=1,\dots,N$. 
\vskip 0.30cm
Undoubtedly, the first point is that two $\omega_{\text{hyp}}$-admissible metrics on
a fixed line bundle differe only by a constant factor.
The second point is that we have already had arithmetically natural admissible
metrics on $K_M$ and ${\Cal O}_M(P_i)$,  i.e., the corresponding
$\omega_{\text{hyp}}$-Arakelov metrics $\rho_{\text{Ar},\omega_{\text{hyp}}}$ and 
$\rho_{\text{Ar},\omega_{\text{hyp}},P_i}$, $i=1,\dots,N$, respectively.
Hence, the geometrically natural admissible metrics on $K_M$ and ${\Cal
O}_M(P_i)$ we seek should be proportional to the corresponding arithmetically
natural admissible metrics defined by using hyperbolic Green's functions.
\vskip 0.30cm
With this in mind, we define the geometrically natural
admissible metric on $K_M$ by multiplying the $\omega_{\text{hyp}}$-Arakelov
metric $\rho_{\text{Ar},\omega_{\text{hyp}}}$  the factor
$A_{\text{Ar,hyp}}(M^0)$. Denote the resulting Hermitian line bundle by
$\underline{K_M}_{\text{hyp}}$. That is, we have
$$\underline{K_M}_{\text{hyp}}:=\underline{K_M}_{\omega_{\text{hyp}}}\cdot
A_{\text{Ar,hyp}}(M^0),\eqno(6.2.1)$$ or equivalently,
$$\rho_{\text{hyp};K_M}:=\rho_{\text{Ar},\omega_{\text{hyp}}}\cdot
A_{\text{Ar,hyp}}(M^0).\eqno(6.2.2)$$
\vskip 0.30cm
Once a geometrically meaningful admissible metric is introduced on $K_M$, we are
left with the only problem to define a similar metric on the cuspidal line bundle.
For this purpose, we introduce the following additional principal:
for our theory of metrics, all the punctures should be viewed as the same, i.e.,
there should be no difference when we impose the geometrically meaningful
admissible metrics $\rho_{\text{hyp};P_i}$ on ${\Cal O}_M(P_i)$, $i=1,\dots,N$ by
modifying $\rho_{\text{Ar},\omega_{\text{hyp}},P_i}$'s. In other words,
from now on,  we assume that the (resulting constant) ratio
$$C^i_{\text{hyp}}:=e^{c^i_{\text{hyp}}}:=
\rho_{\text{hyp};P_i}/\rho_{\text{Ar};\omega_{\text{hyp}};P_i}\eqno(6.2.3)$$
does not depend on $i$. Obviously, with all this, the condition in (6.1.2),
claiming that
$\underline{K_M(P_1+\dots+P_N)}_{\omega_{\text{hyp}}}$ multiplying by
$e^{a_{\text{hyp}}+c^1_{\text{hyp}}+\dots+c^N_{\text{hyp}}}$ is
isometric to 
$K(P_1+\dots+P_N)$
together with the natural metric $\tau^\vee_{\text{hyp};K_M(P_1+\dots+P_N)}$
 induced from $\tau_{\text{hyp}}$ on $M^0$, 
determines  the constant $c_{\text{hyp}}:=c^i_{\text{hyp}}$, $i=1,\dots,N$
and hence the metrics $\rho_{\text{hyp};P_i}$ on ${\Cal O}_M(P_i)$, $i=1,\dots,N$,
uniquely. This then finishes our discussion on how to impose geometrically
meaningful admissible metrics on $K_M$ and on the cuspidal line bundles
respectively. For our own convinences, we set
$$\underline{{\Cal O}_M(P_i)}_{\text{hyp}}:=({\Cal
O}_M(P_i),\rho_{\text{hyp};P_i}),\qquad i=1,\dots,N.\eqno(6.2.4)$$
\vskip 0.30cm
\noindent
{\it Remark 6.2.1.} Here we omit the condition (6.1.3), as 
this can hardly be the case.  Nevertheless, by the above discussion, we 
see that the ratio of the metrics on both hands of (6.1.3) is a constant 
which depends only on punctured Riemann surface $M^0$ itself. 
\vskip 0.30cm
\noindent
(6.3) Before finally giving the geometric interpretation for our metric on
 determinant of cohomology,   we in this
subsection using the result in (5.4)  give the Mumford type isometry associated to
hyperbolic metrics, which will be used in the next section on Takhtajan-Zograf
metrics.
\vskip 0.30cm
For this purpose,  we apply Theorem 5.4.1 as follows. First of all, take $\omega$
to be the normalized hyperbolic volume $\omega_{\text{hyp}}$.
Secondly, set
$(\alpha;\beta_1,\dots,\beta_N)$ in subsection (5.4) to be
$(a_{\text{Ar,hyp}};c^1_{\text{hyp}},\dots,c^N_{\text{hyp}})$ introduced in (6.2).
Finally,  denote the resulting corresponding Hermitian line bundles by the underline
with the lower index hyp, e.g., $\underline{\lambda_n}_{\text{hyp}}$ stands for
$\overline{\lambda_n}^{\alpha;\beta}$,
$\underline{\Delta_n}_{\text{hyp}}$ stands for
$\overline{\Delta_n}^{\alpha;\beta}$ , etc.. Then we have the following
\vskip 0.30cm
\noindent
{\bf Theorem 6.3.1.} (Fundamental Theorem with respect to Hyperbolic Metrics)
{\it
With the same notation as above, for any positive integer $n$, we have the following
isometries:
\vskip 0.30cm
\noindent
(a) (Serre isometry) 
$$\underline{\lambda_n}_{\text{hyp}}\simeq\underline{\lambda_{1-n}}_{\text{hyp}};$$
\vskip 0.30cm
\noindent
(b) (Generalized Mumford isometry) 
$$\underline{\lambda_n}_{\text{hyp}}^{\otimes 12}\simeq
\underline{\Delta_0}_{\text{hyp}}^{\otimes 6n^2-6n+1}\otimes
\underline{\Delta_1}_{\text{hyp}}
\otimes\underline{\Delta_2}_{\text{hyp}}^{\otimes
-12n+10}\otimes{\Cal O}(e^{a(q)});$$
\vskip 0.30cm
\noindent
(c) (Generalized Mumford isometry) 
$$\underline{\lambda_n}_{\text{hyp}}\simeq
\underline{\lambda_1}_{\text{hyp}}^{\otimes 6n^2-6n+1}\otimes
\underline{\Delta_1}_{\text{hyp}}^{\otimes -{{n(n-1)}\over 2}}
\otimes\underline{\Delta_2}_{\text{hyp}}^{\otimes(n-1)^2}\otimes{\Cal
O}(e^{-{{n(n-1)}\over 2}\cdot a(q)}).$$}
\vskip 0.30cm
Obviously, even though we only discuss our metric theory for a single curve, but
the technique can be globalized so that we can apply the above discussion for
a family of curves. In particular, this then works  over the Teichm\"uller space
$T_{q,N}$ of punctured Riemann surfaces of signature $(q,N)$ as well as over the
moduli space ${\Cal M}_{q,N}$ of punctured Riemann surfaces of signature $(q,N)$.
Moreover, as
$$\underline{K_M(P_1+\dots+P_N)}_{\text{hyp}}\simeq
(K_M(D),\tau_{\text{hyp};K_M(D)}^\vee),\eqno(6.3.1)$$
 by a work of Wolpert ([Wo1]), (see e.g. [TZ2] and (7.2)
below for the detail,) we may deduce that
$$c_1(\underline{\Delta_0}_{\text{hyp}})={{\omega_{\text{WP}}}\over{\pi^2}}.
\eqno(6.3.2)$$ Here $\omega_{\text{WP}}$ denotes the Weil-Petersson K\"ahler
form.  Thus in
particular, we have the following
\vskip 0.30cm
\noindent
{\bf Corollary 6.3.2.} {\it With the same notation as above, for all
positive integers $n$, we have the following identities of (1,1)-forms on
$T_{q,N}$ and hence on ${\Cal M}_{q,N}$:
 $$12\,c_1(\underline{\lambda_n}_{\text{hyp}})=(6n^2-6n+1)
{{\omega_{\text{WP}}}\over{\pi^2}}+c_1(
\underline{\Delta_1}_{\text{hyp}})-(12n-10)c_1(
\underline{\Delta_2}_{\text{hyp}}).\eqno(6.3.3)$$}
\vskip 0.20cm
Theorem 6.3.1 and Corollary 6.3.2 will be used to connect our work with the
beautiful pioneer work of Taktajan and Zograf ([TZ1,2]) in \S7. 
\vskip 0.30cm
\noindent
(6.4) The geometric interpretation of our metrics on determinant of
cohomology at this stage is given in terms of the new metric on $\lambda(K_M)$
with respect to the hyperbolic metric. 
\vskip 0.30cm
Realize $M^0$ as a quotient $\Gamma\backslash {\Cal H}$ of the upper
half-plane by the action of a torsion free finitely generated Fuchsian group
$\Gamma$. Then it is well-known that we may choose $\Gamma\subset PSL(2,{\Bbb R})$
to be a subgroup generated by $2q$ hyperbolic transformations
$A_1,B_1,\dots,A_q,B_q$ and
$N$ parabolic transformtions $S_1,\dots,S_N$ satisfying the single relation
$$A_1B_1A_1^{-1}B_1^{-1}\dots A_qB_qA_q^{-1}B_q^{-1}S_1\dots S_N=1.$$ Choose
 a normalized basis of abelian differentials
$\psi_1,\dots,\psi_q$, i.e., a basis of the vector space $H^0(M,K_M)$ so
that 
$$\int_z^{A_iz}\psi_j(w)dw=\delta_{ij},\quad
\int_z^{B_iz}\psi_j(w)dw=:\tau_{ij},\quad i,j=1,\dots,q,$$ with $\delta_{ij}$
the Kronecker symbol and $\tau=(\tau_{ij})$ the period matrix of $M$.
\vskip 0.30cm
On $\lambda(K_M)$, choose the section $(\psi_1\wedge\dots\wedge
\psi_q)\otimes 1^\vee$, with $1$ the canonical section of $H^1(M,K_M)\simeq
{\Bbb C}$. Then we have the following
\vskip 0.30cm
\noindent
{\bf Theorem 6.4.1.} {\it With the same notation as above, 
as the metric on $\lambda(K_M)$,}
$$\eqalign{~&\langle  (\psi_1\wedge\dots\wedge
\psi_q)\otimes 1^\vee,(\psi_1\wedge\dots\wedge
\psi_q)\otimes 1^\vee\rangle_{h_{\underline
{K_M}_{\text{hyp}}}(\underline
{K_M}_{\text{hyp}})}\cr
=&\Big({{\text{det}\,(\text{Im}\tau})\cdot {2\pi(2q-2)}}\Big)\cdot
 (\text{det}^*(\Delta_{\text{hyp}}))^{-1}.\cr}$$
\vskip 0.30cm
\noindent
{\it Proof.} From \S5, we know that the new metrics on determinant of
cohomology obey the rules in \S3 for the Riemann-Roch metrics. 
Thus we see that
$$\eqalign{~&h_{\underline{K_M}_{\text{hyp}}}(\underline
{K_M}_{\text{hyp}})\cr
=&h_{\underline
{K_M}_{\omega_{\text{hyp}}}}(\underline
{K_M}_{\omega_{\text{hyp}}})\cdot e^{{{2q-2}\over {12}}\cdot
a_{\text{Ar,hyp}}}\quad(\text{by}\ (3.2.4))\cr =&h_{\underline
{K_M}_{\omega_{\text{hyp}}}}(\underline
{K_M}_{\omega_{\text{hyp}}})
\cdot e^{{{2q-2}\over {12}}\cdot{{12}\over
2}\cdot {1\over {2q-2}}\cdot\Big(\text{log}
{{\text{det}^*\Delta_{{\text{Ar}}}}\over {A_{\text{Ar}}(M)}}
-\text{log}{{\text{det}^*\Delta_{{\text{hyp}}}}\over
{2\pi(2q-2)}}\Big)}\cr
~&\quad(\text{by}\ (6.1.8))\cr
=&h_{\underline
{K_M}_{\omega_{\text{hyp}}}}(\underline
{K_M}_{\omega_{\text{hyp}}})\cdot \sqrt{
{{\text{det}^*\Delta_{{\text{Ar}}}}\over {A_{\text{Ar}}(M)}}}\cdot
\Big(\sqrt{{{\text{det}^*\Delta_{{\text{hyp}}}}\over
{2\pi(2q-2)}}}\Big)^{-1}\cr
=&\Big(h_{\underline
{K_M}_{\omega_{\text{Ar}}}}(\underline
{K_M}_{\omega_{\text{Ar}}})\cdot \sqrt{
{{\text{det}^*\Delta_{{\text{Ar}}}}\over {A_{\text{Ar}}(M)}}}\Big)\cdot
\Big(\sqrt{{{\text{det}^*\Delta_{{\text{hyp}}}}\over
{2\pi(2q-2)}}}\Big)^{-1}\cr
~&\quad(\text{by}\ (3.3.2)).\cr}$$ But we know that $h_{\underline
{K_M}_{\omega_{\text{Ar}}}}(\underline
{K_M}_{\omega_{\text{Ar}}})$ is simply the Quillen metric on $\lambda(K_M)$
with respect to the Arakelov metric, thus,
by definition, 
$$h_{\underline
{K_M}_{\omega_{\text{Ar}}}}(\underline
{K_M}_{\omega_{\text{Ar}}})=h_{F;\omega_{\text{Ar}}}\cdot\Big(\sqrt{
{{\text{det}^*\Delta_{{\text{Ar}}}}\over
{A_{\text{Ar}}(M)}}}\Big)^{-1}.$$ This to say, 
$$h_{\underline
{K_M}_{\omega_{\text{Ar}}}}(\underline
{K_M}_{\omega_{\text{Ar}}})\cdot \sqrt{
{{\text{det}^*\Delta_{{\text{Ar}}}}\over {A_{\text{Ar}}(M)}}}$$
is simply the Faltings metric $h_{F;\omega_{\text{Ar}}}$ on $\lambda(K_M)$,
which is nothing but the determinant of the $L^2$-pairing on $H^0(M,K_M)$.
Therefore, by Serre duality, we see that
$$\eqalign{~&\langle (\psi_1\wedge\dots\wedge
\psi_q)\otimes 1^\vee,(\psi_1\wedge\dots\wedge
\psi_q)\otimes
1^\vee\rangle_{h_{\underline
{K_M}_{\text{hyp}}}(\underline
{K_M}_{\text{hyp}})}\cr
=&\langle (\psi_1\wedge\dots\wedge
\psi_q)\otimes 1^\vee,(\psi_1\wedge\dots\wedge
\psi_q)\otimes
1^\vee\rangle_{h_{F;\omega_{\text{Ar}}}}\cdot 
\Big(\sqrt{{{\text{det}^*\Delta_{{\text{hyp}}}}\over
{2\pi(2q-2)}}}\Big)^{-2}\cr
=&\Big({\text{det}\,\text{Im}\tau}\cdot {2\pi(2q-2)}\Big)\cdot
 (\text{det}^*(\Delta_{\text{hyp}}))^{-1}.\cr}$$ This completes the proof of
the theorem.
\vskip 0.30cm
As a direct consequence of the proof of Theorem 6.4.1, we have the following
\vskip 0.30cm
\noindent
{\bf Corollary 6.4.2.} {\it With the same notation as above, if $M^0=M$, i.e.,
if $M^0$ is compact, then $\underline{K_M}_{\text{hyp}}$ is nothing but $K_M$
together with standard (smooth) hyperbolic metric. In other words, when the Riemann
surface is compact, then the Arakelov-Poincar\'e volume is simply the ratio between
the standard (smooth) hyperbolic metric and the Arakelov metric with respect to the
normalized hyperbolic volume form.}
\vskip 0.30cm
\noindent
{\it Remark 6.4.1.} Recall that by various Mean Value Lemmas, for our matric
theory, the $\omega$-Arakelov metric is essentially the original Arakelov
metric, which is in the nature of Euclidean geometry. Hence, the above corollary
shows that the Arakelov-Poincar\'e volume  indeed measures how far the Euclidean
aspect of a compact Riemann surface is away from its Poincar\'e aspect.
\vskip 0.30cm
\noindent
{\we \S7. On Takhtajan-Zograf metric over moduli space of punctured Riemann
surfaces} 
\vskip 0.30cm
\noindent
(7.1) To faciliate ensuing discussion on an application of our metric, we
in this subsection recall some results of Takhtajan and Zograf ([TZ1,2]). 
\vskip 0.30cm
For a punctured Riemann surface $M^0$ of signature $(q,N)$ (with $2q+N\geq 3$),
 let $\Gamma$ be a torsion free Fuchsian group unformizing $M^0$, i.e., 
$M^0\simeq \Gamma\backslash {\Cal H}$, where ${\Cal H}$ denotes the complex upper-half plane.
Denote by $\Gamma_1,\dots,\Gamma_N$ the set of non-conjugate parabolic subgroups in $\Gamma$, and
for every $i=1,\dots,N$, fix an element $\sigma_i\in PSL(2,{\Bbb R})$ such that 
$\sigma_i^{-1}\Gamma_i\sigma_i=\Gamma_\infty$, where the group $\Gamma_\infty$ is 
generated by the parabolic transformation $z\mapsto z+1$.
As usual, define the Eisenstein series $E_i(s,z)$ corresponding to the $i$-th cusp of the group 
$\Gamma$ for $\text{Re}(s)>1$ by
$$E_i(s,z):=\sum_{\gamma\in\Gamma_i\backslash \Gamma}\text{Im}(\sigma_i^{-1}\gamma z)^s,
\quad i=1,\dots,N.\eqno(7.1.1)$$
\vskip 0.20cm
Denote the Teichm\"uller space of punctured Riemann surfaces of signature $(q,N)$ by $T_{q,N}$.
Then at the point $[M^0]$ corresponding to a punctured Riemann surface $M^0$, the tangent space 
$T_{[M^0]}T_{q,N}$ can be naturally identified with the space  $\Omega^{-1,1}(M^0)$ of harmonic $L^2$-tensors on $M^0$ of type
(-1,1). Define the Weil-Petersson metric on $T_{q,N}$ by
$$\langle \psi,\psi\rangle_{\text{WP}}:=
\int_{M^0}\psi\bar\psi d\mu_{\text{hyp}},\eqno(7.1.2)$$ where
$\psi,\psi\in 
\Omega^{-1,1}(M^0)$ are considered as tangent vectors of $T_{q,N}$ at $[M^0]$, and 
$d\mu_{\text{hyp}}=2\pi(2q-2+N)\,\omega_{\text{hyp}}$ is the K\"ahler
form  corresponding to the uniformizing hyperbolic metric
$\tau_{\text{hyp}}$ induced from ${\Cal H}\to\Gamma\backslash  {\Cal
H}\simeq M^0$ with Gaussian curvature -1.
\vskip 0.30cm
Following Takhtajan and Zograf, for $i=1,\dots,N$, define the metric
$\langle\,,\,\rangle_i$ on
$T_{q,N}$ by setting $$\langle \phi,\psi\rangle_i:=
\int_{M^0}\phi\bar\psi E_i(\cdot,2) d\mu_{\text{hyp}},
\quad\phi,\,\psi\in\Omega^{-1,1}(M^0).\eqno(7.1.3)$$
Here $E_i(z,s)$ is the Eisenstein series defined in (7.1.1). In [TZ2], it is
proved that $\langle\,,\,\rangle_i$,
$i=1,\dots,N$, are K\"ahler metrics on $T_{q,N}$. Moreover, their sum
$\sum_{i=1}^N
\langle\,,\,\rangle_i$ gives a new K\"ahler metric $\langle\,,\,\rangle_{\text{cusp}}$, the cusp metric or the 
Takhtajan-Zograf metric, on $T_{q,N}$, which is invariant under the action of the Teichm\"uller modular group.
Denote the corresponding K\"ahler form by $\omega_{\text{cusp}}$, or
$\omega_{\text{TZ}}$.
\vskip 0.30cm
For compact Riemann surfaces $M$, a work of D'Hoker-Phong [D'HP] and
Sarnak [Sa] shows that  the regularized determinant
$\text{det}^*\Delta_n$ associated to
$K_M^{\otimes n}$ with  respect to hyperbolic metrics defined via the zeta
function formalism of Ray-Singer, is equal,  up to a constant  multiplier
depending only on $q$ and $n$, to $Z_M'(1)$ for $n=0,1$,  and $Z_M(n)$ for
$n\geq 2$. Here $Z_M(s)$ denotes the Selberg zeta function associated to
$M$. Motivated by this and the Quillen metric on determinant of cohomology, for
punctured Riemann surfaces, Takhtajan and Zograf ([TZ1,2])  define
$\text{det}_{\text{TZ}}^*\Delta_n$ with respect to hyperbolic metrics by
simply setting
$$\text{det}_{\text{TZ}}^*\Delta_n:=\cases Z_{M^0}'(1),&$if\ $n=0,1$$;\\
Z_{M^0}(n),&$if\ $n\geq 2$$.\endcases$$  Here $Z_{M^0}(s)$
denotes the Selberg zeta function of $M^0$ defined in (6.1.4). Moreover, for any
$n\in {\Bbb Z}_{\geq 1}$,  on
$\lambda_n:=\lambda(K_{M}^{\otimes n}\otimes{\Cal O}_{M}(P_1+\dots+P_N)^{\otimes
(n-1)})$, they introduce the norm
$h_{\text{TZ},n}$ by setting $$h_{\text{TZ},n}:=h_P\cdot
\text{det}_{\text{TZ}}^*\Delta_n^{-{1\over 2}},$$ where $h_P$ denotes the
determinant  of Petersson norm on $\lambda_n$. 
\vskip 0.30cm
\noindent
{\bf Theorem 7.1.1.} ([TZ1,2]) {\it With the same notation as
above, as (1,1) forms on $T_{q,N}$:
$$c_1(\lambda_n,h_{\text{TZ},n})={{6n^2-6n+1}\over{12}}
{{\omega_{\text{WP}}}\over{\pi^2}}-{1\over
9}\omega_{\text{TZ}}.\eqno(7.1.4)$$}
\vskip 0.30cm
Note that $Z_{M^0}(s)$ has a simple zero at $s=1$, 
our definition (6.1.7) of the regularized determinant of
Laplacian $\Delta_{\text{hyp}}$ differs from  Takhtajan
and Zograf's $\text{det}_{\text{TZ}}^*\Delta_1$ up to a universal constant factor
depending only on the signature $(q,N)$. (See e.g., Remark 6.1.2.)  Therefore, by
Theorem 6.4.1, we have the following
\vskip 0.30cm
\noindent
{\bf Theorem 7.1.1$'$.}  {\it With the same notation as
above, as (1,1) forms on $T_{q,N}$:
$$12\cdot c_1(\underline{\lambda_1}_{\text{hyp}})=
{{\omega_{\text{WP}}}\over{\pi^2}}-{4\over
3}\omega_{\text{TZ}}.\eqno(7.1.5)$$}
\vskip 0.20cm
Moreover, as all bundles, forms, and
metrics are invariant under the action of the Teichm\"uller modular group, 
(7.1.4) and (7.1.5) actually induce the same relations on the moduli space ${\Cal
M}_{q,N}$ of punctured Riemann surfaces of signature $(q,N)$ (in the sense of
$V$-manifolds).
\vskip 0.30cm
\noindent
(7.2) Now let us look at Theorem 6.3.1(b) and (c) carefully.
First of all, by definition in (6.1), we know that
$\overline{K_M}^{a_{\text{hyp}}}(\overline
{P_1}^{c_{\text{hyp}}}+\dots+\overline{P_N}^{c_{\text{hyp}}})$
is really nothing but $K_M(P_1+\dots+P_N)$ together with the hyperbolic metric
$\tau^\vee_{\text{hyp};K_M(P_1+\dots+P_N)}$
naturally induced from the hyperbolic volume form $d\mu_{\text{hyp}}$.
Thus in particular, over
$T_{q,N}$, the Teichm\"uller space of punctured Riemann surfaces of
signature $(q,N)$, the first Chern form for the  metrized Deligne pairing 
$$\langle \overline{K_M}^{a_{\text{hyp}}}(\overline
{P_1}^{c_{\text{hyp}}}+\dots+\overline{P_N}^{c_{\text{hyp}}}),
\overline{K_M}^{a_{\text{hyp}}}(\overline
{P_1}^{c_{\text{hyp}}}+\dots+\overline{P_N}^{c_{\text{hyp}}})\rangle$$
may be naturally associated to ${{\omega_{\text{WP}}}\over {\pi^2}}$ by the work of
Wolpert as stated in [TZ2]. In fact, we have the following
\vskip 0.30cm
\noindent
{\bf Theorem 7.2.1.} ([Wo], [TZ2])  {\it With the
same notation as above,  
 $$\int
c_1(K_M(P_1+\dots+P_N),\tau^\vee_{\text{hyp};K_M(P_1+\dots+P_N)})^2={{\omega_{\text{WP}}}\over{\pi^2}}.$$}
\vskip 0.30cm
On the other hand, by arithmetic intersection theory, we have
$$c_1(\langle \bar L,\bar L'\rangle)=\int c_1(\bar L)\cdot c_1(\bar L').$$ (See
e.g., [El], where the metrics involved are supposed to be smooth. But this
restriction can be easily removed   to apply here, as the
singularities of our metrics on admissible metrized line bundles are of
hyperbolic growth.) Hence, 
$$c_1\Big(\langle \overline{K_M}^{a_{\text{hyp}}}(\overline
{P_1}^{c_{\text{hyp}}}+\dots+\overline{P_N}^{c_{\text{hyp}}}),
\overline{K_M}^{a_{\text{hyp}}}(\overline
{P_1}^{c_{\text{hyp}}}+\dots+\overline{P_N}^{c_{\text{hyp}}})\rangle\Big)
={{\omega_{\text{WP}}}\over {\pi^2}}.$$
As a direct consequence, as claimed in Corollary 6.3.2,
$$12c_1(\underline{\lambda_n}_{\text{hyp}})=(6n^2-6n+1){{\omega_{\text{WP}}}\over
{\pi^2}}+c_1(\underline{\Delta_1}_{\text{hyp}})
-(12n-10)c_1(\underline{\Delta_2}_{\text{hyp}}).
\eqno(7.2.1)$$
Moreover, bu the generalized Mumford isometry,
$$\Big(\underline{\lambda_n}_{\text{hyp}}\otimes \underline{\Delta_2}_{\text{hyp}}^{\otimes
n-1}\Big)^{\otimes
12}\simeq
\underline{\Delta_0}_{\text{hyp}}^{\otimes
6n^2-6n+1}\otimes
\underline{\Delta_1}_{\text{hyp}}\otimes
\underline{\Delta_2}_{\text{hyp}}^{\otimes
-2}\otimes {\Cal
O}(e^{a(q)})\quad\text{for\ all}\ n\geq 1.\eqno(7.2.2)$$ 
Here
$$\underline{\Delta_0}_{\text{hyp}}:=\langle\underline{K_N({\Bbb
P}_1+\dots+{\Bbb P}_N)}_{\text{hyp}},
\underline{K_N({\Bbb
P}_1+\dots+{\Bbb
P}_N)}_{\text{hyp}}\rangle(\pi_N).$$ 
\vskip 0.30cm
\noindent
(7.3) Put  (7.2.2) in the language of differential forms, we see
that, if $n\geq 1$,
$$12c_1\Big(\underline{\lambda_n}_{\text{hyp}}\otimes
\underline{\Delta_2}_{\text{hyp}}^{\otimes
n-1}\Big)=(6n^2-6n+1)c_1\Big(\underline{\Delta_0}_{\text{hyp}}\Big)
+c_1\Big(\underline{\Delta_1}_{\text{hyp}}\otimes
\underline{\Delta_2}_{\text{hyp}}^{\otimes
-2}\Big).\eqno(7.3.1)$$ Let $n=1$,
and use Theorem 7.1.2 and Theorem 7.2.1, we have the following
\vskip 0.30cm
\noindent
{\bf Theorem 7.3.1.} (Fujiki-Weng) {\it With the same notation
as above,
 $${4\over
3}\omega_{\text{TZ}}=c_1\Big(\underline{\Delta_1}_{\text{hyp}}\otimes
\underline{\Delta_2}_{\text{hyp}}^{\otimes
-2}\Big).$$ So ${4\over
3}\omega_{\text{TZ}}$ can be realized as the first Chern form of the metrized line
bundle $\underline{\Delta_1}_{\text{hyp}}\otimes
\underline{\Delta_2}_{\text{hyp}}^{\otimes
-2}$. In particular, the Takhtajan-Zograf metric is K\"ahler and
${4\over 3}\omega_{\text{TZ}}$ is a Hodge metric form.}
\vskip 0.30cm
\noindent
{\it Remark 7.3.1.} Note that as a line bundle $\Delta_2$ is indeed  trivial. So if
we only interested in the results in the bundle version, we in fact have the
following simple relation
$$\lambda_n^{\otimes 12}\simeq \Delta_0^{\otimes
6n^2-6n+1}\otimes\Delta_1.\eqno(7.3.2)$$ On the other hand, arithmetically,
$\underline{\Delta_2}_{\text{hyp}}$ is far from being trivial. It seems to be
equally interesting to study the associated smooth function on ${\Cal M}_{q,N}$
resulting from the corresponding metric on
$\underline{\Delta_2}_{\text{hyp}}$.
 \vskip 0.30cm
\noindent
{\bf Theorem 7.3.2.} (Fundamental Relations for Riemann Surfaces) {\it With the same
notation as above, ${\Cal M}_{q,N}$, for $n\geq 0$, there exists the canonical
isometry, up to a constant factor depending only $q,N$ and $n$,
$$(\lambda_n,h_{\text{TZ},n})^{\otimes
12}\simeq 
\underline{\Delta_0}_{\text{hyp}}^{\otimes 6n^2-6n+1}\otimes\Big(
\underline{\Delta_1}_{\text{hyp}}\otimes
\underline{\Delta_2}_{\text{hyp}}^{\otimes
-2}\Big).\eqno(7.3.3)$$ In particular, up to a
constant factor depending only $q,N$ and $n$,
$$(\lambda_n,h_{\text{TZ},n})=\underline{\lambda_n}_{\text{hyp}}\otimes
\underline{\Delta_2}_{\text{hyp}}^{\otimes
n-1}.\eqno(7.3.4)$$}
\vskip 0.30cm
\noindent
{\it Proof.} Easily we see that by (7.3.1) and Theorem 7.3.1,  both sides
of (7.3.3) have the same first Chern form over ${\Cal M}_{q,N}$. Thus by the
structure of the stably compactified moduli space $\overline{{\Cal M}_{q,N}}$, 
the metrics on both sides are proportional to each other. (See e.g. the proof of
Theorem A2.4.1 in the appendix.) This completes the proof of the theorem. 
\vskip 0.30cm
\noindent
\centerline {\we APPENDIX:}
\vskip 0.30cm
\centerline{\we  Arithmetic Factorization Theorem in terms of
Intersection}
\vskip 0.30cm
In this appendix, we  propose an arithmetic factorization for 
Weil-Petersson  geometry, Takhtajan-Zograf  geometry and Selberg geometry
associated to punctured Riemann surfaces. Unlike the rest of this
paper, the discussion here is rather informal, in particular, not so many rigorous
proofs are given for  the assertions. So for the time being, the reader may
simply understand them as some working hypothesis. On the other hand, we
anticipate that this arithmetic factorization will play a key role in studying
the global geometry of moduli spaces of Riemann surfaces.
\vskip 0.30cm
\noindent
{\we \S A1.  Degeneration of Weil-Petersson metrics}
\vskip 0.30cm
\noindent
(A1.1) We start with Masur's result on degeneration of Weil-Petersson
metrics. Let
${\Cal M}_q$ be the moduli space of compact Riemann surfaces of genus $q\geq 2$.
Denote its stably compactification by $\overline{{\Cal M}_q}$. Let
$p\in E:=\overline{{\Cal M}_q}\backslash {\Cal M}_q$ be a boundary point and let
$U$ be a small neighborhood of $p$. Let $\pi:\Delta^n\to\Delta^n/G=U$ be a local
uniformizing chart with holomorphic coordinates
$((z_i)_{i=1}^r,(w_j)_{j=1}^{3q-3-r})$ such that $\pi((0),(0))=p,$ and 
$\pi^{-1}(U\cap E)=\cup_iz_i^{-1}(0).$ As before, denote by $\omega_{\text{WP}}$
the Weil-Petersson K\"ahler form. Write
$$\eqalign{~&\pi^*(\omega_{\text{WP}}|_{U\cap {\Cal M}_q})\cr
=&\sqrt{-1}\sum a_{i\bar j}dz_i\wedge d\bar z_j-2\text{Im}\sum b_{i\bar k}
dz_i\wedge d\bar w_k+\sqrt{-1}\sum c_{k\bar l} dw_k\wedge d\bar w_l.\cr}$$ Then we
have the following fundamental result of Masur ([Ma]):
\vskip 0.30cm
\noindent
{\bf Theorem A1.1.1.} (Masur) {\it For $((z_i),(w_k))$ near $0$,
\vskip 0.30cm
\noindent
(i) $\displaystyle{{{C^{-1}}\over {|z_i|^2(-\log|z_i|)^3}}\leq a_{i\bar i}\leq 
{{C^{-1}}\over {|z_i|^2(-\log|z_i|)^3}}}$, for a constant $C>0$;
\vskip 0.30cm
\noindent
(ii) $\displaystyle{a_{i\bar
j}=O\Big({1\over{|z_i|\,|z_j|\,(-\log|z_i|)^3(-\log|z_j|)^3}}\Big)}$, if
$i\not=j$;
\vskip 0.30cm
\noindent
(iii) $\displaystyle{b_{i\bar
k}=O\Big({1\over{|z_i|(-\log|z_i|)^3}}\Big)}$;
\vskip 0.30cm
\noindent
(iv) $\displaystyle{\lim_{((z_i),(w_k))\to 0}c_{k\bar l}=h_{k\bar l}}$ for a
certain constant positive definite hermitian matrix $(h_{k\bar l})$.}
\vskip 0.30cm
\noindent
(A1.2) Previously, when  mathematicians talked about Masur's result, they
usually paid much more attention on the first three conclusions, i.e., the
asymptotic behavior. Even though such an asymptotic behavior is very important
geometrically, but it has little arithmetic meanning.
On the other hand, there is another very important
part which is hiden in conclusion (iv) -- this last statement clearly indicate 
that the restriction of the Weil-Petersson K\"ahler form to the boundary could
result a metric on the boundary. Indeed, later we will see that it is not
unreasonable to expect that this induced form should coincide with the
Weil-Petersson K\"ahler form for the boundary.
\vskip 0.30cm
\noindent
(A1.3) To justify what we said in (A1.2), we next recall yet another fundamental
result due to Wolpert ([Wo2]).
\vskip 0.30cm
Let $\overline{\pi_q}:\overline{{\Cal C}_q}\to \overline{{\Cal M}_q}$ (resp.
${\pi_q}:{{\Cal C}_q}\to {{\Cal M}_q}$) be the
universal curve over the stably moduli space (resp. moduli space) of compact
Riemann surfaces. Then it is well-known that the standard hyperbolic matric on
compact Riemann surfaces may be glued together to give a smooth metric
on the relative canonical line bundle $K_{\pi_q}$. On the other hand, even on
singular fibers of $\overline{\pi_q}$, we may get standard hyperbolic metrics on
the corresponding punctured Riemann surfaces. A natural question is whether these
(singular) hyperbolic metrics can also be glued together so that we can get a
certain type of metric on the relative canonical line bundle
$K_{\overline{\pi_q}}$. The answer is yes. In fact, we have the following
\vskip 0.30cm
\noindent
{\bf Theorem A1.3.1.} (Wolpert) {\it With the same notation as above, the
resulting metric on the relative canonical line bundle $K_{\overline{\pi_q}}$
obtaining from standard hyperbolic metrics on the fibers is continuous and good.}
\vskip 0.30cm
As a direct consequence of this result, Wolpert then deduces that,
in the sense of currents, on the compactified moduli space $\overline {{\Cal
M}_q}$,
${{\omega_{\text{WP}}}\over {\pi^2}}$ is the curvature form of a continuous
metric $h_{\text {WP}}$ on a certain line bundle and the metric $h_{\text{WP}}$
may be approximated by smooth positive curvature metrics. We later will give an
alternative proof of this  statement.
\vskip 0.30cm
\noindent
(A1.4) Motivated by Wolpert's result, in [TW2], we study
how to glue admissible metrics along with a degeneration family of compact Riemann
surfaces. This roughly goes as follows.
\vskip 0.30cm
To facilitate ensuing discussion, we first recall the plumbing construction of
a degenerating family of Riemann surfaces starting from $M$ as follows (cf.
e.g. [Fay1], [Ma] and [Wo2]). Let $M^0:=M\backslash\{p\}$. Then $M^0$ is a punctured
Riemann surface with two punctures $p_1,\,p_2$ in place of $p$, where
$p_1,\,p_2$ correspond to two points in the  normalization $\tilde M$ of
$M$. Denote the unit disc in ${\Bbb C}$ by $\Delta$. For $i=1,\,2$, fix a
coordinate function $z_i:U_i\to \Delta$ such that $z_i(p_i)=0,$ where $U_i$
is an open neighborhood of $p_i$. For each $t\in \Delta$, let $S_t:=\{(x,y)\in
\Delta^2:xy=t\}$. Now for each $t\in\Delta$, remove the discs
$|z_i|<|t|,\,i=1,\,2$, from $M$ and glue the remaining surface with $S_t$ via
the identification
$$z_1\sim (z_1,t/z_1)\ \ \text{and}\ \ z_2\sim(t/z_2,z_2).\eqno(A1.4.1)$$
The resulting surfaces $\{M_t\}_{t\in \Delta}$ form an analytic family
$\pi:{\Cal M}\to\Delta$ with $M_0=M$. Here $\pi$ denotes the holomorphic
projection map. Note that for $t\not=0$, each fiber $M_t$ is a compact
Riemann surface of genus $q$. Also the node $p$ does not disconnect the
Riemann surface when removed from $M$. The restriction of $\text{
ker}(d\pi)$ to ${\Cal M}\backslash \{p\}$ forms a holomorphic line bundle
over ${\Cal M}\backslash \{p\}$ such that $L|_{M_t}=TM_t$ and
$L|_{M^0}=TM^0$, which will be called the vertical line bundle. Note that
$\text{ker}(d\pi)$ itself does not form a line bundle over ${\Cal M}$ since
$\text{ker}(d\pi)$ is of rank 2 at $p$. Similarly, one may construct
a degenerating family of compact Riemann surfaces such that the center fiber
is a nodal curve with a separating node.
\vskip 0.45cm
Now we are ready to state the following result of To and myself in [TW2].
\vskip 0.30cm
\noindent
Let $\{M_t\}$ be a family of compact
Riemann surface of genus $q\geq 2$ degenerating to a Riemann surface $M$
of genus $q-1$ with a single  node $p$ as described above. Let $L=\{L_t\}$
be a line bundle on $\{M_t\}$. Then
\vskip 0.30cm
\noindent
(i) in the case when $M_0$ is with a non-separating node, there is a continuous
metric $\rho$ defined everywhere on $\{M_t\}$, except possibly at the node, such
that 
\vskip 0.30cm
\noindent
(a) the restriction of $\rho$ to $\{M_t\}_{t\not=0}$ is smooth; 
\vskip 0.30cm
\noindent
(b) for each
$t\not=0$, the restriction of $\rho$ to $L_t$ is $d\mu_{\text{hyp},t}$-admissible;
\vskip 0.30cm
\noindent
(ii) in the case when $M_0$ is with a separating node, the following two
conditions are equivalent: 
\vskip 0.30cm
\noindent
(A) there is a
continuous metric $\rho$ defined everywhere on $\{M_t\}$, except possibly at the
node, such that 
\vskip 0.30cm
\noindent
(a) the restriction of $\rho$ to $\{M_t\}_{t\not=0}$ is smooth; 
\vskip 0.30cm
\noindent
(b) for each
$t\not=0$, the restriction of $\rho$ to $L_t$ is $d\mu_{\text{hyp},t}$-admissible;
\vskip 0.30cm
\noindent
(B) the degrees $d_1$ and $d_2$ of $L_0$ on the irrediucible components
$M_0^{(1)}$ of genus $q_1$ and $M_0^{(2)}$ of genus $q_2$ of $M_0$ satisfy
$$d_1(2q_2-1)=d_2(2q_1-1).$$
\vskip 0.30cm
\noindent
Moreover we know that, in general, if we do not have the degree condition as in
(B), admissible metrics on
$\{L_t\}_{t\not=0}$ will continuously extend to one of the irreducible
component, while blow up to infinity on the other irreducible component.
\vskip 0.30cm
\noindent
{\we \S A2.  Arithmetric Factorization Theorem: a proposal}
\vskip 0.30cm
\noindent
(A2.1) We start with an algebraic factorization theorem.
So we go back to study the universal curve
$\overline{\pi_{q,N}}:\overline{{\Cal C}_{q,N}}\to \overline{{\Cal M}_{q,N}}.$
On $\overline{{\Cal C}_{q,N}}$, the following line bundles are well-defined:
$K$, the relative canonical line bundle; ${\Cal O}({\Bbb P}_i)$, $i=1,\dots,N$,
$N$ sections. Hence, by using Deligne pairing formalism, we get the following
line bundles over  $\overline{{\Cal M}_{q,N}}$:
$\langle K({\Bbb P}_1+\dots+{\Bbb P}_N),K({\Bbb P}_1+\dots+{\Bbb P}_N)\rangle$;
$\langle K({\Bbb P}_1+\dots+{\Bbb P}_N),{\Bbb P}_i\rangle$;
$\langle K,{\Bbb P}_i\rangle$, $i=1,\dots,N$. Parallelly, we have Mumford type
line bundles $\lambda_n$ introduced in \S5. Moreover, we know that these
line bundles on $\overline{{\Cal M}_{q,N}}$ satisfy Momford type relations, i.e.,
Theorem 5.2.2 (on ${\Cal M}_{q,N}$).
\vskip 0.30cm
Now we want to know how these bundles change when we restrict them to the
boundary of $\overline{{\Cal M}_{q,N}}$, or better when we  pull back  these
bundles via the normalization of the stable curves. For simplicity, we only study
the case when  one more non-separating node is involved. So we have the following
natural map
$\alpha:{{\Cal M}_{q-1,N+2}}\rightarrow\overline{{\Cal
M}_{q,N}}$. We will use $\widetilde K$ to denote the relative canonical line bundle
for the universal curve on ${\Cal M}_{q-1,N+2}$, and use
$\widetilde {\Bbb P}_i$, $i=1,\dots,N$ and
${\Bbb R},\,{\Bbb S}$ to denote $N+2$ sections (so that ${\Bbb P}_i$ corresponds
to $\widetilde {\Bbb P}_i$, $i=1,\dots,N$ and ${\Bbb R},\,{\Bbb S}$ are two more
sections corresponding to the non-separating node for the restriction of
the original universal curve $\overline{\pi_{q,N}}$ to the boundary.)
\vskip 0.30cm
Obviously, in this case, we have the following algebraic factorization:
\vskip 0.30cm
\noindent
(a) the line bundle $\langle K({\Bbb P}_1+\dots+{\Bbb P}_N),K({\Bbb
P}_1+\dots+{\Bbb P}_N)\rangle$ changes to 
$\langle \widetilde K(\widetilde {\Bbb P}_1+\dots+\widetilde {\Bbb P}_N+{\Bbb R}+{\Bbb
S}),\widetilde K(\widetilde{\Bbb P}_1+\dots+\widetilde{\Bbb P}_N+{\Bbb R}+{\Bbb S})\rangle$;
\vskip 0.30cm
\noindent
(b) the line bundle $\langle K({\Bbb P}_1+\dots+{\Bbb P}_N),{\Bbb P}_i\rangle$
changes to
$\langle \widetilde K(\widetilde {\Bbb P}_1+\dots+\widetilde {\Bbb P}_N+{\Bbb R}+{\Bbb
S}),\widetilde{\Bbb P}_i\rangle$, $i=1,\dots,N$;
\vskip 0.30cm
\noindent
(c) the line bundle
$\langle K,{\Bbb P}_i\rangle$ changes to 
$\langle \widetilde K,\widetilde{\Bbb P}_i\rangle\otimes \langle{\Bbb R}+{\Bbb
S},\widetilde{\Bbb P}_i\rangle$, $i=1,\dots,N$;
\vskip 0.30cm
\noindent
(d) the line bundle $\lambda_n^{\otimes 12}\otimes (\langle K({\Bbb
P}_1+\dots+{\Bbb P}_N),{\Bbb P}_1+\dots+{\Bbb P}_N\rangle)^{\otimes 6(n-1)}$
changes to
$\widetilde\lambda_n^{\otimes 12}\otimes \Big(\langle \widetilde K(\widetilde{\Bbb
P}_1+\dots+\widetilde{\Bbb P}_N+{\Bbb R}+{\Bbb S}),\widetilde{\Bbb
P}_1+\dots+\widetilde{\Bbb P}_N+{\Bbb R}+{\Bbb S}\rangle\Big)^{\otimes
6(n-1)}\otimes \Big(\langle \widetilde K,{\Bbb R}+{\Bbb
S}\rangle\otimes\langle\widetilde K({\Bbb R}+{\Bbb
S}),{\Bbb R}+{\Bbb
S}\rangle\Big).$
\vskip 0.30cm
The proof may be obtained by looking at the intersection first, which gives (a),
(b) and (c). With (a), (b) and (c), (d) is a direct consequence of the generalized
Mumford relation from Theorem 5.2.2. In fact note that many components in (d) are
trivial line bundles, we may rewrite (d) as
\vskip 0.30cm
\noindent
(d$'$) the line bundle $\lambda_n^{\otimes 12}$
changes to
$\widetilde\lambda_n^{\otimes 12}\otimes \langle \widetilde K,{\Bbb R}+{\Bbb
S}\rangle.$
\vskip 0.30cm
\noindent
(A2.2) With the above discussion, we may now offer the following 
 global picture for the geometry of punctured Riemann surfaces, or
better, for the geometry of moduli spaces of punctured Riemann surfaces.
\vskip 0.30cm
First of all, our generalized Mumford type isometrie in Theorem 6.3.1, together with
Theorem 7.3.1 and Theorem 7.3.2
expose explicitly the intrinsic relations among the spectrum geometry given by
Selberg zeta functions, the deformation geometry given by Weil-Petersson metric,
and the cusp geometry given by Eisenstin series via Takhtajan-Zograf metrics.
\vskip 0.30cm
Secondly, the deformation geometry and the cusp geometry are in the nature of
arithmetic intersection theory. So their properties should be relatively easier to
understand. As a consequence, via  Mumford isometries, we then could get
information about the spectrum geometry, which is in nature of cohomology theory.
\vskip 0.30cm
Finally, algebraic factorization (a) shows that via the degeneration
and normalization process, the Weil-Petersson geometry factors extremely well. So
the arithmetic counter part should be established in a rather formal way.
Similarly, we can apply this comment to the cusp geometry by looking at algebraic
factorizations (c) and (d).
\vskip 0.30cm
\noindent
(A2.3) Now we indicate how one can do the arithmetic factorization for deformation
geometry, i.e., the Weil-Petersson metric. To this end, we need to get a simliar
result as in Theorem A1.3.1. More precisely, the following statement should be first
established.
\vskip 0.30cm
\noindent
{\it With the same notation as above, the
resulting metric on the relative logarithmic canonical line bundle
for $\overline{\pi_{q,N}}$ obtaining from the standard (singular) hyperbolic
metrics on the fibers is continuous and good.}
\vskip 0.30cm
Indeed, as pointed out by Professor Fujiki, this result may be obtained as a
direct consequence of Theorem A1.3.1 of Wolpert by using the structure of the
universal curve over $\overline{{\Cal M}_{q,N}}$.
\vskip 0.30cm
With this, in particular, we will obtain a continuous metric on the
Deligne pairing $\langle K({\Bbb P}_1+\dots+{\Bbb P}_N),K({\Bbb
P}_1+\dots+{\Bbb P}_N)\rangle$. That is, we have a continuous metrized line bundle
$\langle \underline{K({\Bbb P}_1+\dots+{\Bbb
P}_N)}_{\text{hyp}},\underline{K({\Bbb P}_1+\dots+{\Bbb
P}_N)}_{\text{hyp}}\rangle$ on $\overline{{\Cal M}_{q,N}}$. From here, by using
the property of Deligne metric on Deligne pairings, we can further conclude that in
Masur's result, e.g., Theorem A1.1.1(iv), the positive metrix is really nothing but
the one coming from the Weil-Petersson metric for
${\Cal M}_{q-1,N+2}.$ So the Weil-Petersson metric factors extremely well.
\vskip 0.30cm
\noindent
(A2.4) To give the arithmetic factorization for Takhtajan-Zograf metric, we need
first decomposite them into $N$-pieces. In fact, by symmetry, we see that
$$c_1\Big(\langle \underline{K}_{\text{hyp}},\underline{{\Bbb
P}_i}_{\text{hyp}}\rangle\otimes\langle \underline{K({\Bbb P}_1+\dots+{\Bbb
P}_N)}_{\text{hyp}}, \underline{{\Bbb P}_i}_{\text{hyp}}\rangle\Big)={4\over
3}\omega_{\text{WZ}}^{(i)},
\quad i=1,\dots,N.$$ Here, for  $i=1,\dots,N$, $\omega_{\text{WZ}}^{(i)}$ denotes
the $i$-th Takhtajan-Zograf K\"ahler form associated to the $i$-th
Takhtajan-Zograf  metric on ${\Cal M}_{q,N}$ defined by using the $i$-th
Eisenstein series. Thus we should consider line bundles
$\langle K,{\Bbb P}_i\rangle\otimes\langle K({\Bbb P}_1+\dots+{\Bbb P}_N), {\Bbb
P}_i\rangle$, $i=1,\dots,N$. By algebraic factorization (b) and (c) in (A2.1), we
see that they factor into $\Big(\langle
\widetilde K,\widetilde {\Bbb P}_i\rangle\otimes\langle \widetilde
K(\widetilde{\Bbb P}_1+\dots+\widetilde{\Bbb P}_N+{\Bbb R}+{\Bbb S}),
\widetilde {\Bbb P}_i\rangle\Big)\otimes\Big(\langle{\Bbb R}+{\Bbb S},
\widetilde {\Bbb P}_i\rangle\Big)$,
$i=1,\dots,N$. At this moment, I should say that arithmetically, the appearance
of $\langle{\Bbb R}+{\Bbb S}, \widetilde {\Bbb P}_i\rangle$ is extremely
unpleasent, as I cannot show that arithmetically it is trivial. (Indeed, I think
it is hardly the case.) But the line bundle $\langle{\Bbb R}+{\Bbb S},
\widetilde {\Bbb P}_i\rangle$ is trivial, so let us for the time being pretend that
such an appearance is harmless for the discussion follows.
\vskip 0.30cm
As for the case about the Weil-Petersson metric, to understand the arithmetic
factorization, we need to study the corresponding metrics on line bundles over
$\overline{{\Cal C}_{q,N}}$ first. In (6.1), we already introduce natural
metrics on $K_M$ and on $P_i$, $i=1,\dots,N$ for each punctured Riemann surface by
introducing an invariant called Arakelov-Poincar\'e volume. The point now is whether
such  metrics will form  continuous metrics on $K$ and ${\Bbb P}_i$,
$i=1,\dots,N$, when we are working on a family. By looking at the result of To and
myself recalled in (A1.4), it is resonable to conclude that if the total volume of
$M$ with respect to the metric induced from  $\underline{K_M}_{\text{hyp}}$ is an
absolute constant, i.e., the associated total volume is independent of
$M$ and $P_1,\dots,P_N$, we should then can glue these metrics  together to get 
continuous metrics on $K$ and ${\Bbb P}_i$, $i=1,\dots,N$ on $\overline{{\Cal
C}_{q,N}}$. For this latest purpose, we next give the  final main result of
this paper. 
\vskip 0.30cm
\noindent
{\bf Theorem A2.4.1.} {\it With the same notation as above, the total volume of
$M$ for the metric induced from $\underline{K_M}_{\text{hyp}}$ is a constant
depending only on $q$ and $N$.}
\vskip 0.30cm
\noindent
{\it Proof.} This is a direct consequence of the geometric interpretation of our
determinant metric on $\underline{\lambda_1}_{\text{hyp}}$. Indeed, if we denote
the total volume of $M$ for the metric induced from
$\underline{K_M}_{\text{hyp}}$  by
$A(M;\underline{K_M}_{\text{hyp}})$, then up to a constant
depending only on $q$ and $N$,  the inner product of our determinant metric for
the generator $1\otimes(\phi_1\wedge\dots\wedge \phi_q)^\vee$ of $\lambda_1$ is
nothing but 
$A(M;\underline{K_M}_{\text{hyp}})$ times the
inverse $\text{det}^*\Delta_{\text{hyp}}$. Here, as before, we denote
$\{\phi_i\}_{i=1}^q$ an orthonormal basis of $H^0(K_M)$ with respect to the
natural pairing.
 But we know that
$\text{det}^*\Delta_{\text{hyp}}$ up to a constant
depending only on $q$ and $N$ is simply $Z'(1)$ with $Z(s)$ denoting the
corresponding Selberg zeta function, by the fact that $Z(s)$ has a simple zero at
$s=1$. Thus by the curvature formula for $\underline{\lambda_1}_{\text{hyp}}$, we
see that
$$dd^cA(M;\underline{K_M}_{\text{hyp}})=0,$$ when we move $M$ in ${\Cal M}_{g,N}$.
But in $\overline{{\Cal M}_{q,N}}$,
locally, the absolute horizontal boundary and the relative horizontal boundaries,
i.e., the fake diagonal divisors $\Delta_S$ can be contracted. (I learn this from
Prof. Fujiki on the way to Kinosaki.) This completes the proof.
\vskip 0.30cm
\noindent
{\it Remark A2.4.1.} Indeed, we would like to guess that
$A(M;\underline{K_M}_{\text{hyp}})=2\pi(2q-2)$. But for the time being, it seems
to be quite imporssible to prove this, as we need more precise degeneration
information for the quantities introduced in this paper. On the other hand, if
this is true, then there is a great chance to simplify the discussion in \S6 and
\S7. 
\vskip 0.30cm
\noindent
{\it Remark A2.4.2.} The reader should know that the metric defined on
$\underline{K_M}_{\text{hyp}}$  is obtained in an arithmetic manner: we first
use the hyperbolic Green's function and the associated beta function to define an
hyperbolic Arakelov metric on $K_M$, which is quite suitable for our arithmetic
purpose; then we multiple this metric by a highly transcendental invariant,
the so-called Arakelov-Poincar\'e volume to obtain the metric, which is motivated
by the work of D'Hoker-Phong and Sarnak. So the possible geometric definition of
this  metric,  proposed in the previous remark, on $K_M$ would be very interesting.
\vskip 0.30cm
With the above, we may assume that there are globally defined metrized line
bundles
$\underline{K}_{\text{hyp}}$ and $\underline{{\Bbb P}_i}_{\text{hyp}}$,
$i=1,\dots,N$, with  continuous metrics on $\overline {{\Cal C}_{q,N}}$. Hence we
further get metrized line bundles $\langle
\underline{K}_{\text{hyp}},\underline{{\Bbb
P}_i}_{\text{hyp}}\rangle\otimes\langle \underline{K({\Bbb P}_1+\dots+{\Bbb
P}_N)}_{\text{hyp}},\underline{{\Bbb P}_i}_{\text{hyp}}\rangle$,
$i=1,\dots,N$ with continuous metrics on $\overline {{\Cal M}_{q,N}}$. This 
clearly shows that the Takhtajan-Zograf metrics
$\langle\cdot,\cdot\rangle_{\text{TZ}}^{(i)}$, $i=1,\dots,N$, have natural
factorizations, which is exactly the same as what has happened  for Weil-Petersson
metrics. In particular, we see that the degereration of the $i$-th Eisenstein
series will correspond to exactly the $i$-th Einstein series of the central fiber.
We would like to point out that such a degeneration has been studied by 
others, notably, Wolpert ([Wo3]). But the picture drew here by using the
arithmetic factorization seems to be  more clear then what is obtained before. In
fact, we see that the additional Eisenstein series corresponding to new
punctures $R$ and $S$ can never  be the limit of the original Eisenstein
series from  nearby fibers. Nevertheless, by the arithmetic
factorization for the determinant line bundles
$\underline{\lambda_n}_{\text{hyp}}$, which may now be obtained by using the
above arithmetic factorization of arithmetic intersection via the generalized
Mumford type isometries, we see that the additional Eisenstein series
corresponding to new punctures $R$ and $S$ are obtained from the spectrum
degeneration via  Selberg zeta functions. So inseatd of traditionally studying
the degeneration of the combination of the Selberg zeta function and the small
eigen-values, one may directly study the degeneration of the Selberg zeta
function itself, as we expect this will give additional information for
additional Eisenstein series corresponding to new punctures $R$ and $S$.
\vskip 0.30cm
\noindent
(A2.5) We conclude this appendix and hence this paper by the following remark. In
the paper, we offer a way to understand the global geometry of a general Riemann
surface. Undoubted, this is just the beginning of the story.  Personally, I believe
that with the arithmetic factorization proposed here
one may finally find an alternative way to understand the arithmetic Miyaoka-Yau
inequality, if we take Belyi's theorem [Be] and Hilbert's irreducibility theorem
[La1] into consideration, instead of trying to establish a
$p$-adic deformation theory, if it exists.
\vfill
\eject
\noindent
\centerline {\we REFERENCES}
\vskip 0.30cm
\item {[Ai]}  W. Aitken: {\it An arithmetic Riemann-Roch theorem for singular
arithmetic surfaces}, Memoirs of the AMS, No. {\bf 573}, 1996
\vskip 0.30cm
\item {[Ar]} S. Arakelov: Intersection theory of divisors on an arithmetic
surface, Izv. Akad.  Nauk SSSR Ser. Mat., 38 No. {\bf 6}, (1974)
\vskip 0.30cm
\item {[BM]} A. Beilinson and V. Ginsburg,
Infinitesimal structure of moduli spaces of
$G$-bundles, IMRN, Duke J. Math. {\bf 4} (1992),
63-74\vskip 0.30cm
\item {[BM]} A. Beilinson and Y. Manin, The Mumford form and the Polyakov
measure in string theory, Comm. Math. Phys, {\bf 107}, 359-376, (1986)
\vskip 0.30cm
\item{[Be]} G. Belyi, On Galois extensions of a maximal cyclotomic field, Math.
USSR Izv. {\bf 14} (1980), 247-256
\vskip 0.30cm
\item {[De1]}  P. Deligne: {\it Equations diff\'erentielles a points singuliers
r\'eguliers}, Lecture Notes in Math. {\bf 163}, Berlin-Heideberg-New York,
Springer, (1970)
\vskip 0.30cm
\item {[De2]}  P. Deligne: Le d\'eterminant de la cohomologie, {\it
Current trends in arithmetic algebraic geometry},
Contemporary Math. Vol. {\bf 67}, 93-178, (1987) 
\vskip 0.30cm
\item {[DM]}  P. Deligne and D. Mumford, The
irreducibility of the space of curves of given
genus, IHES Publ. Math. {\bf 36} (1969), 75-109
\vskip 0.30cm 
\item{[D'HP]} E. D'Hoker and D.H. Phong, On determinant of Laplacians on
Riemann surfaces, Comm. Math. Physics, {\bf 104}, 537-545, (1986)  
\vskip 0.30cm 
\item{[Ef]} I. Efrat, determinant of Laplacians on surfaces of finite
volume, Commun. Math. Phys, {\bf 119}, 443-451, (1988)
\vskip 0.30cm 
\item{[El]} R. Elkik, M\'etriques sur les fibr\'es d'intersection, Duke Math. J.
{\bf 61}, 303-328, (1990) 
\vskip 0.30cm
\item {[Fa]} G. Faltings: Calculus on Arithmetic Surfaces, Ann. Math. 
{\bf 119}, 387-424, (1984)
\vskip 0.30cm
\item {[Fay1]} J. Fay, {\it Theta functions on Riemann
surfaces}, Lecture notes in Math. {\bf 352}.  
Berlin-Heideberg-New York, Springer, (1973)
\vskip 0.30cm
\item {[Fay2]} J. Fay, {\it Kernel functions, analytic torsion, and moduli
spaces}, Memoirs of AMS, No. {\bf 464}, 1992
\vskip 0.30cm
\item {[Fu]} A. Fujiki, An $L^2$ Dolbeault lemma and its applications, Publ.
RIMS, {\bf 28}, 845-884 (1992)
\vskip 0.30cm
\item {[FM]} W. Fulton and R. MacPherson, A
compactification of configuration spaces, Ann.
Math. (2), {\bf 139} (1994), no. 1, 183-225 
\vskip 0.30cm
\item {[Jo]} J. Jorgenson, Asymptotic behavior of
Faltings delta function. Duke Math. J. {\bf 61}
(1990), 221-254
\vskip 0.30cm
\item {[JL1]} J. Jorgenson and R. Lundilius, Convergence theorems for
relative spectral functions on hyperbolic Riemann surfaces with
hyperbolic cusps, Duke Math.
J. {\bf 80}, 785-819, (1995)
\vskip 0.30cm 
\item{[JL2]} J. Jorgenson and R. Lundelius, Continuity of
relative hyperbolic spectral theory through metric degeneration, Duke Math.
J. {\bf 84}, 47-81, (1997)
\vskip 0.30cm
\item {[Kn]} F. Knudsen,  The projectivity of the
moduli space of stable curves, II, and III, Math.
Scand. {\bf 52}(1983), 161-199, 200-212
\vskip 0.30cm
\item {[KM]} F. Knudsen and D. Mumford,  The
projectivity of the moduli space of stable
curves, I. Math. Scand. {\bf 39} (1976), 19-55
\vskip 0.30cm
\item{[La1]} S. Lang, {\it Fundamentals of Diophantine Geometry}, Springer-Verlag,
Berlin-Heide-berg-New York, 1983
\vskip 0.30cm
\item{[La2]} S. Lang,  $\,${\it Introduction to Arakelov theory},
$\,$Springer-Verlag,
$\,$Berlin-Heideberg-New York, 1988
\vskip 0.30cm
\item{[Ma]} H. Masur, The extension of the
Weil-Petersson metric to the boundary of
Teichm\"uller space, Duke Math. J., {\bf 43}
(1976), 623-635
\vskip 0.30cm 
\item {[MB]} L. Moret-Bailly: Le formulae de Noether pour les
surfaces arithm\'etiques, Invert. Math., (1990)
\vskip 0.30cm 
\item {[Mu1]} D. Mumford, Hirzebruch's proportionality theorem in the
non-compact case, Invent. Math. (1977)
\vskip 0.30cm 
\item {[Mu2]} D. Mumford, Stability of projective
varieties, L'Ens. Math., {\bf 24} (1977), 39-110
\vskip 0.30cm 
\item {[Mu3]} D. Mumford, Towards an Enumerative
geometry of the moduli space of curves, in {\it
Arithmetic geometry}, dedicated to Shafarevich on
his 60th birthday,  271-328, (1983)
\vskip 0.30cm
\item {[Qu]} D. Quillen, determinant of Cauchy-Riemann operators over
Riemann surfaces, Func. Anal. Appl. {\bf 19}, 31-34 (1986)
\vskip 0.30cm
\item{[RS]} D. Ray and I. Singer, Analytic torsion for complex manifolds,
Ann. Math. {\bf 98}, 154-177 (1973) 
\vskip 0.30cm
\item{[Sa]} P. Sanark, determinant of Laplacians, Commun. Math. Phys. {\bf
110}, 113-120 (1987)
\vskip 0.30cm
\item{[So]} Ch. Soul\'e, Géométrie d'Arakelov des surfaces
arithmétiques. Séminaire Bourbaki, Vol.
1988/89. Astérisque No. {\bf 177-178} (1989), Exp.
No. 713, 327--343
\vskip 0.30cm
\item{[TW1]} W.-K. To and L. Weng, The asymptotic behavior of Green's
functions for quasi-hyperbo -lic metrics on degenerating Riemann surfaces,
Manuscripta Math. {\bf 93}, 465-480 (1997)
\vskip 0.30cm
\item{[TW2]} W.-K. To and L. Weng,  Admissible Hermitian metrics on families of line bundles over
degenerating Riemann surfaces, preprint, 1998
\vskip 0.30cm
\item{[TZ1]} L. Takhtajan and P. Zograf, The Selberg zeta function and a new
K\"ahler metric on the moduli space of punctured Riemann surfaces, J. Geo.
Phys. {\bf 5}, 551-570 (1988)
\vskip 0.30cm
\item{[TZ2]} L. Takhtajan and P. Zograf, A local index theorem for
families of $\bar\partial$-operators on punctured Riemann surfaces and a new
K\"ahler metric on their moduli spaces, Commun. Math. Phys. {\bf
137}, 399-426 (1991)
\vskip 0.30cm
\item{[We1]} L. Weng, $\Omega$-admissible theory, to appear in the Proceedings of
London Math. Soc.
\vskip 0.30cm
\item{[We2]} L. Weng, {\it Hyperbolic Metrics,
Selberg Zeta Functions and Arakelov Theory for Punctured
Riemann Surfaces,} Lecture Notes
at Osaka University, manuscript, Osaka University (1998)
\vskip 0.30cm
\item{[Wo1]} S. Wolpert, Chern forms and the Riemann tensor for the moduli
space of curves, Inv. Math. {\bf 85}, 119-145 (1986)
\vskip 0.30cm
\item{[Wo2]} S. Wolpert, The hyperbolic metric
and the geometry of universal curve. J. Differ.
Gemo. {\bf 31} (1990), 417-472
\vskip 0.30cm
\item{[Wo3]} S. Wolpert, Spectral limits for
hyperbolic surfaces, I, II, Invent. Math. {\bf
108} (1992), 67-90, 91-129
\end